\pdfoutput=1
\documentclass{article}
\usepackage[margin=25mm]{geometry}
\usepackage{amsmath}
\usepackage{amsfonts}
\usepackage{amssymb}
\usepackage{graphicx}
\usepackage[utf8]{inputenc}
\usepackage[T2A]{fontenc}
\usepackage[ukrainian,russian,english]{babel}
\usepackage[pdftex]{hyperref}
\usepackage{xcolor}
\usepackage[all]{xy}

\newtheorem{theorem}{Theorem}

\title{Visualization of Morse flow  with two saddles on 3-sphere 
diagrams}
\author{Svitlana Bilun and Alexandr Prishlyak}

\begin{document}

\maketitle
\begin{abstract}
   We describe all possible topological structures of Morse-Smale flows without closed trajectories on a three-dimensional sphere, which have two sources, two sinks, one saddle of Morse index 1, one saddle of Morse index 2, and no more than 10 saddle connections. To classify such flows, a generalized Heegaard diagram or Pr-diagram is used, which in this case consists of a sphere and two closed curves, the intersection points of which correspond to saddle connections. We have found all possible, up to homeomorphism, ways to embed two circles in a 2-sphere with no more than 10 points of transversal intersection and construct its planar visualisations.
\end{abstract}
\section*{Introduction}

In the  paper, a topological classification of Morse flows on the three-dimensional sphere $S^3$ is carried out using diagrams (configurations) consisting of 2 closed curves without self-intersections, which intersect transversally on the sphere in $R^3$. These diagrams are generalized Heegaard diagrams for closed 3-manifolds and Pr-diagrams for 3-manifolds with the boundary. We consider them as a graph of valence 4 with two selected cycles. Since a Morse flow given a consistent ordering of fixed points defines the structure of a Morse function, the given studies can be used to classify these functions. At the same time, the Morse flow has the structure of a gradient flow of the Morse function for some Riemannian metric \cite{Smale1961}.

The use of graph theory to classify dynamical systems is typical in two and three dimensions. For Morse flows on surfaces, the first invariants were called distinguishing graphs \cite{Andronov1937, Leontovich1955, Smale1960, Peixoto1959, Peixoto1973, Palis1982, Palis1970, Fleitas1975}. After \cite{Oshemkov1998}, three-color graphs became popular \cite{prishlyak2020three, ABP2022}. Then these invariants were generalized to a wider class of dynamical systems, and other invariants were also constructed \cite{prislyak2017morse, Prishlyak2019, prishlyak2003topological, prishlyak2003sum,    kkp2013,  prish2002vek,     Prishlyak2017, Prishlyak2022, Prishlyak2021,  Prishlyak2020, Kybalko2018,  Prishlyak2019,   Palis1968,  Giryk1996,   Bolsinov2004, Kadubovskyj2005, Poltavec1995}.

For simple Morse functions on closed surfaces, the main invariant is the Reeb graph\cite{Reeb1946, Kronrod1950}. These graphs are also generalized to a wider class of functions \cite{prishlyak2001conjugacy, hladysh2019simple, hladysh2017topology, prishlyak2002topological1, prishlyak2002morse, prishlyak2000conjugacy,  prishlyak2007classification, lychak2009morse, prish2002Morse, prish2015top, prish1998sopr, bilun2002closed, Bilun13, Sharko1993}.

In dimension 3, diagrams consisting of a surface and the curves embedded in them are often used to classify flows and functions \cite{prish1998vek, prishlyak1999equivalence, prish2001top, prishlyak2002morse1,  Prishlyak2002, prishlyak2003regular, prishlyak2002topological, prishlyak2005complete, Prishlyak2007, Hatamian2020, BPP2022}. For polar Morse flows, they coincide with Heegaard diagrams.

 The main invariants of graphs and their embeddings in surfaces can be founded in \cite{prishlyak1997graphs, Harary69, pontr86, tatt88, HW68, GT87}.

The first section gives the basic definitions and topological properties associated with the diagrams (configurations) of two closed curves on the sphere. For each such diagram, an adjacency graph and its summary matrix are considered, on the basis of which the classification of diagrams (configurations) is carried out. The number of non-equivalent configurations with 2, 4, 6, 8, 10 intersection points is counted, which are equal to 1, 1, 2, 4, 13, respectively.

In the final section, it is shown that the number of topologically non-equivalent Morse flow on the 3-dimensional sphere, which have one fixed point of index 1 and 2 and two sinks and sources each, is equal to the corresponding number of non-equivalent diagrams. In this way, the number of topologically non-equivalent Morse flow is calculated. At the same time, each intersection point on the diagram corresponds to a single integral trajectory connecting special points of index 1 and 2.

We express our gratitude to the students Eugenie Buket, Alex Danciger, Christian Hatamian and Julia Pilina for their essential assistance in performing the calculations and to the president of AUI Stepan Moskalyuk for supporting these studies.

\section{Basic constructions}

We consider the classification of diagrams ( configurations) of two closed curves without self-intersections on the sphere in $R^ 3$ that intersect transversally in the section. Two configurations are considered equal or equivalent if there is a homeomorphism of the sphere onto itself, which translates the curves of one configuration into the curve of the other. For each configuration, a scan can be constructed - a stereographic projection of the sphere onto a plane from an arbitrary point of the sphere that does not belong to any of the curves. It is easy to show that from the equivalence of sweeps follows the equivalence of configurations, but not vice versa. An arc of a curve is a part of a curve bounded by two points of intersection with another curve. We  call a region a part of the surface of the sphere, which is bounded by arcs of curves and to which other regions do not belong. Each configuration can have at most different sweeps than regions on the sphere.

Each configuration of two closed curves on the sphere is matched by an adjacency graph constructed as follows: each region corresponds to a vertex of the graph, the arc of the curve between two regions corresponds to an edge incident to both vertices of the regions. 

Example.

\begin{figure}[ht]  \center{\includegraphics[height=4.5cm]{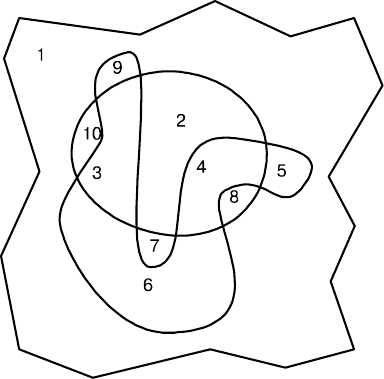}  \ \ \ \ \ \ 
\includegraphics[height=4.5cm]{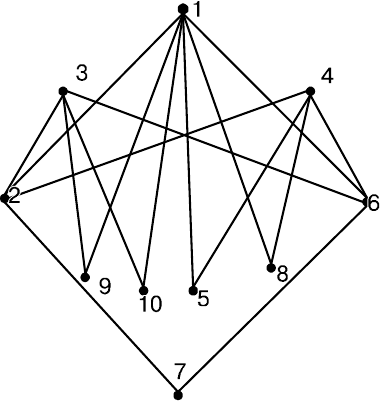}}
\caption{ sweep and the correspondent graph}
 \end{figure}

Regions of the sphere can always be colored in 2 colors so that each arc of the curve separates regions of different colors. Indeed, first we draw the first curve on the sphere and color the two resulting areas. Next, we draw the second curve and repaint the areas bounded by it.
Therefore, we can decolorize the sweep pattern into 2 colors so that each circle arc separates areas of different colors.
It follows that the adjacency graph is bichromatic and bipartite according to the terminology \cite{Harary69, tatt88, HW68, GT87}, i.e. the set of graph vertices can be divided into two subsets (destinies) in such a way that the ends of each edge belong to different destinies. Each edge of the graph (circle arc) can also be assigned responsibility for one of the two
colors according to whether the arc belongs to one or another circle.
Therefore, we will introduce an additional marking of the graph:

\begin{figure}[ht]  \center{\includegraphics[height=4.5cm]{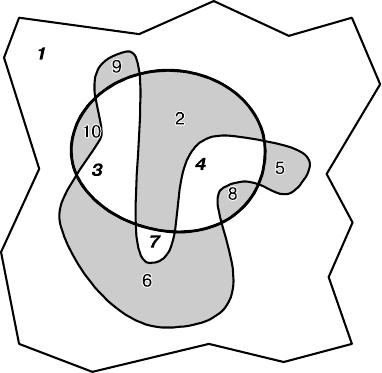}  \ 
\ \ \ \ \ \ \ \includegraphics[height=4.5cm]{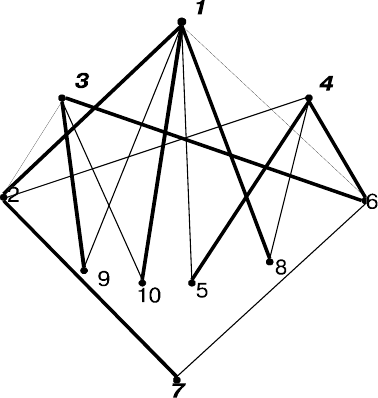}}
 \label{p34} 
 \caption{ marked sweep and the correspondent marked graph }
 \end{figure}

It follows from the construction that the degree of each vertex is even and the same number of edges of different colors are incident to it - the main invariant of the considered graphs.

Each bipartite graph corresponds to and is given by the composite adjacency matrix of the bipartite graph \cite{Harary69}, the matrix of graph or just a matrix. Let $x_1, \dots, x_n$ be the vertices of one color (of one fate), and $y_1, \dots, y_m$ be the vertices of another color (of the second fate). Then in the matrix $M=M_{i,j}$: 

$M_{i,j}=1$ if $x_i$ and $y_j$ is connected by an edge of first color,

$M_{i,j}=-1$ if $x_i$ and $y_j$ is connected by an edge of second color,

$M_{i,j}=-1$ if $x_i$ and $y_j$ is not connected by an edge.

Regarding matrices, the main invariant looks like this: the number of rows and columns of the matrix is even, the sum of the elements of each row and column is zero.

The matrix of the graph that was given in the example:

\includegraphics[height=3.5cm]{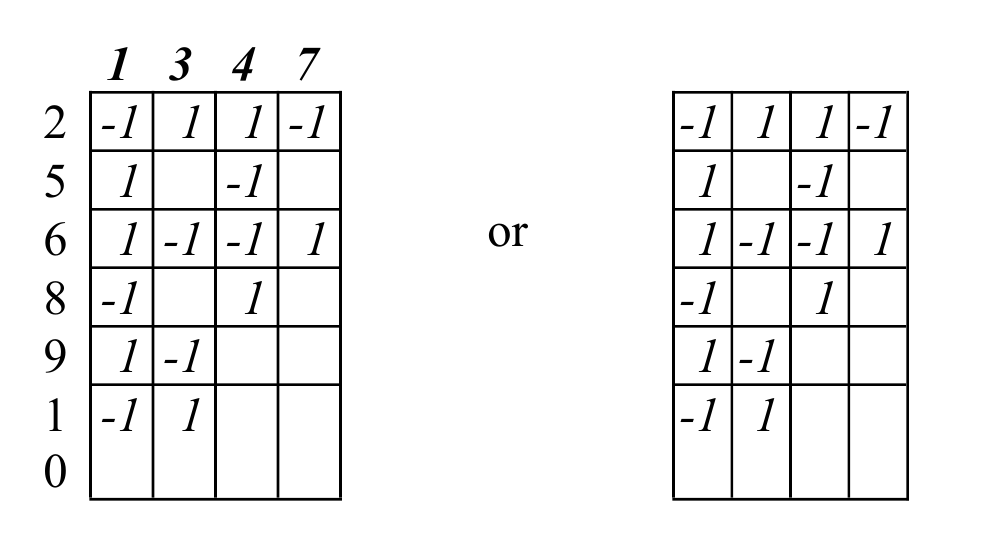}

In our case, the graph is determined by its matrix with accuracy up to arbitrary permutation of rows and columns or rotations by 90 o (corresponds to recoloring of areas) or replacement sign of all variables (corresponds to the recoloring of closed curves). Also, to simplify identification, let's introduce the defining vectors of the graph matrix - two degree vectors of vertices of different colors. For the graph considered above $$(4,4,2,2,2), (6,4,4,2) \ \text{ or } \ \ (4^2,2^3), (6, 4^2,2).$$

The equivalence of two configurations implies the equivalence of their graphs, which is equivalent to the equivalence of the matrices, and the equality of the corresponding determinant vectors follows from it. And vice versa: from the inequality of defining vectors follows the non-equivalence of matrices and graphs, and from it follows the inequality of configurations. We will use determinant vectors to differentiate matrices.

Now let's consider a single crossing operation, with the help of which there is a transition from 2 n points of intersection to 2( n+1 ) points of intersection of curves on the sphere .

Consider an arbitrary region \textbf{1} and two adjacent regions 2 and 3, such that edges \textbf{1} 2 and \textbf{1} 3 are of different colors (belong to different curves). The operation of crossing edges \textbf{1} 2 and \textbf{1} 3 in region \textbf{1} is called intersection 2 1 3. It follows from the construction that operations 2 \textbf{1} 3 and 3 \textbf{1} 2 are the same.

\begin{figure}[ht] 
\center{\includegraphics[height=4.5cm]{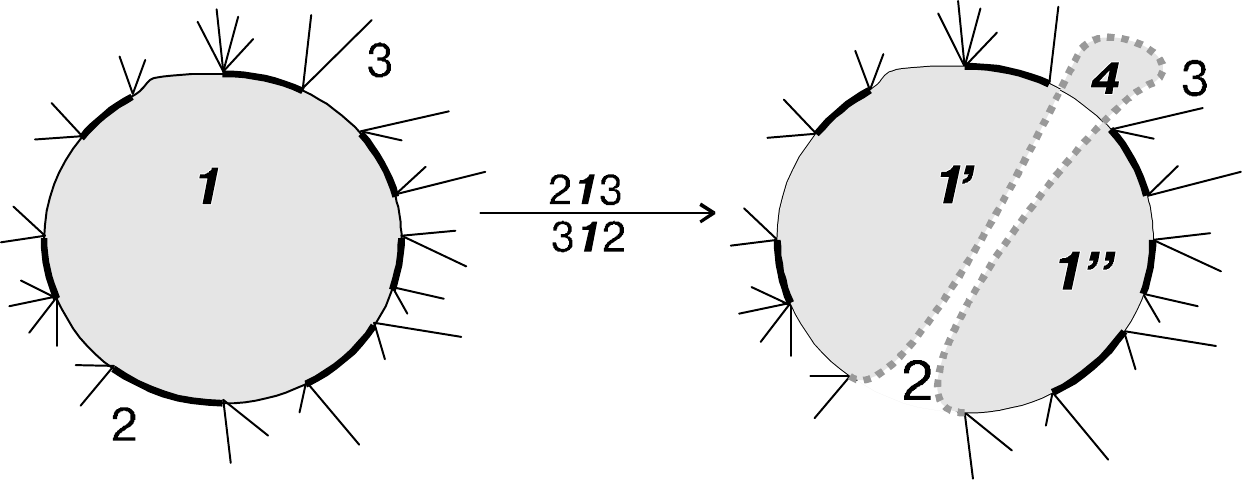}}
\caption{ }
\end{figure}

When crossing, a new area 4 (from two edges) is formed , which is adjacent to areas 3 and 2 and areas 1' and 1'' (obtained as a result of the splitting of region 1 ), which are also adjacent to regions 3 and 2 each. Other areas that were adjacent to area 1 are redistributed (hereinafter this redistribution 
will be called ? - distribution) between areas 1' and 1''. At the same time, areas 1' , 1'' , 4 retain the color of area 1 .

The following happens with the edges: instead of edges 1 2 and 1 3 , edges 1' 2, 1'' 2, 1' 3, 1'' 3 appear , and edges 2 4 and 3 4 also appear , while the color 1 2 remains at edges 1' 2, 1'' 2, 3 4 , and color 1 3 - at edges 1' 3, 1'' 3, 2 4 . Therefore, we have: when intersecting, the number of edges increases by 4, and the number of regions - by 2.

Let's consider what happens to the adjacency graph when intersecting:

\begin{figure}[ht] 
\center{\includegraphics[height=3.5cm]{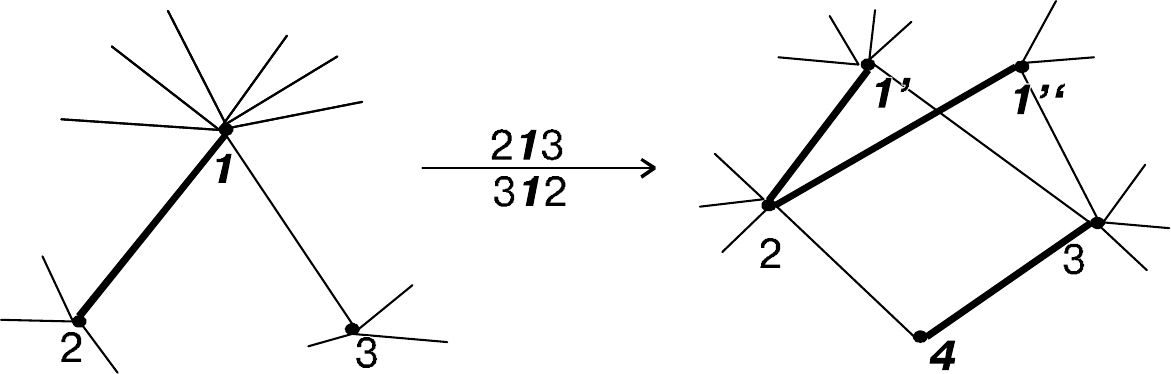}}
\caption{ }
\end{figure}

Vertex 4 is added and instead of vertex 1 two vertices 1' and 1'' are added (while all three new vertices are of the same color) and two pairs of edges of different colors are added. The invariant of parity of degrees of all vertices and the same size of edges of different colors that are incident to one vertex is naturally preserved during interpolation.

Now let's consider the operation of interpolation on the graph matrix: 

In any column (row), we will choose two arbitrary non-zero elements of different signs $x_i$ and $x_j$  in rows (columns) $i$ and $j$, respectively. Let's cross out this column (row). Instead of it, we add 3 columns (rows) to the matrix as follows: in one of them, there are zeros in all places except $i$ and $j$, in the i-th place there is $x_j = - x_i$, in the j-th place there is $x_i = - x_j$; in the other two columns (rows) on the i-th and j-th places are $x_i$ respectively and $x_ j$, and other elements of the crossed out column (row) are distributed ( x -distribution ) in these two columns (rows) in their places according to the invariant.

This matrix shows the 2 1 3 opposition discussed above.

\begin{figure}[ht] 
\center{\includegraphics[height=5.5cm]{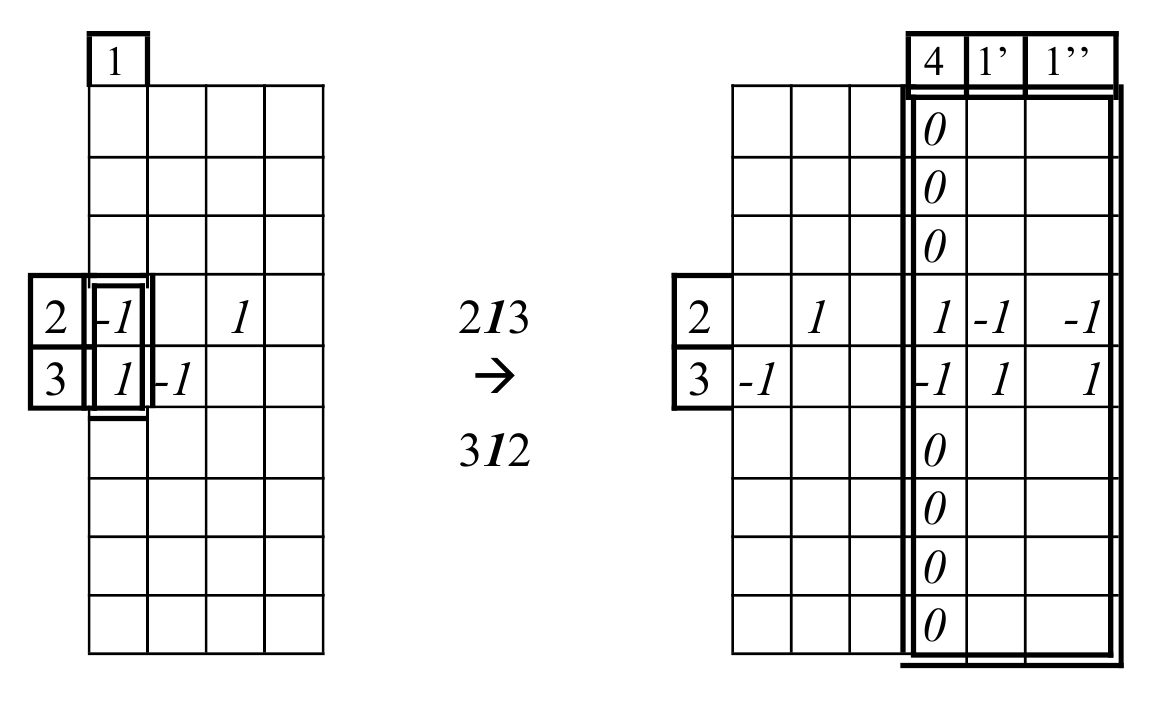}}
\end{figure}

Next, consider the construction of intersections of closed curves on the sphere for 2, 4, 6, 8, 10 intersection points. The number of configurations will be denoted by N(2),..., N(10), respectively.

\newpage 

\section{     Two intersection points}

For this case, we have a single configuration of closed curves on the sphere

\begin{figure}[ht] 
\center{\includegraphics[height=3.5cm]{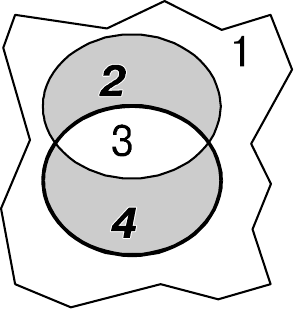}}
\caption{ }
\end{figure}

The corresponding matrix has the form : 

Matrix 2-1

\begin{figure}[ht] 
\center{\includegraphics[height=1.5cm]{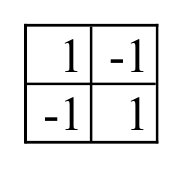}}
\end{figure}

Thus, the following theorem holds.
\begin{theorem}
There is a unique configuration with 2 intersections: $N( 2 )=1.$
\end{theorem}

Next, we will stick to the following algorithm for building configurations:

    1.  for each constructed matrix, we will consider all possible contradictions;
    
    2.  for each extension, we consider all cases of the x-distribution using a sweep-distribution;
    
    3.  based on the new matrix, we construct a new matrix, a new graph and defining vectors, and with the help of the previous sweep we construct a new sweep;
    
    4. we use defining vectors and sweeps to determine the equivalent configuration.

    \section{       Fore intersection points}

From the appearance of the previous matrix, the unity of the matrix, which will be constructed from the previous one by means of multiplication, immediately follows:

\includegraphics[height=1.5cm]{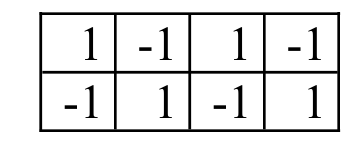}

In this case, the x-distribution is uniform.
Two views of the corresponding graph : 
\begin{figure}[ht] 
\center{\includegraphics[height=3.5cm]{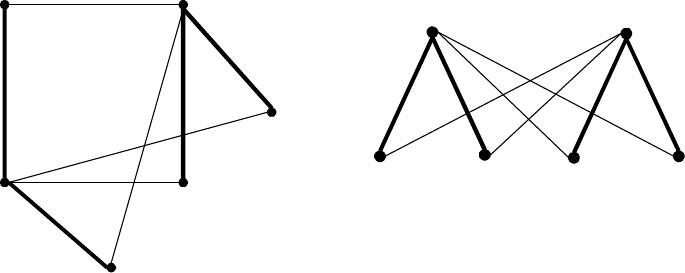}}
\caption{ }
\end{figure}

\hfill
\newpage

This graph corresponds to a single configuration of the intersection of closed curves on the sphere and the following two equivalent sweeps.

\begin{figure}[ht] 
\center{\includegraphics[height=3.5cm]{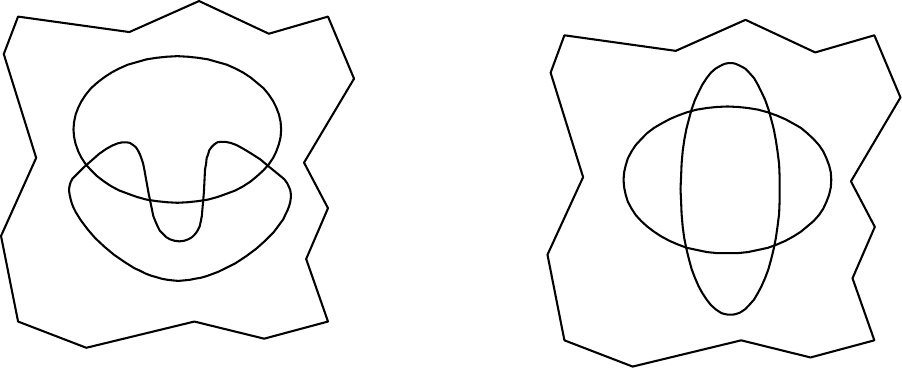}}
\caption{ }
\end{figure}

Thus, the following theorem holds.
\begin{theorem}
There is a unique configuration with 4 intersections: $N( 4 )=1.$
\end{theorem}

          \section{ Six intersection points}
          
In matrices and 4-1 it is possible to construct, without limiting the generality, two intersections.

\includegraphics[height=3.5cm]{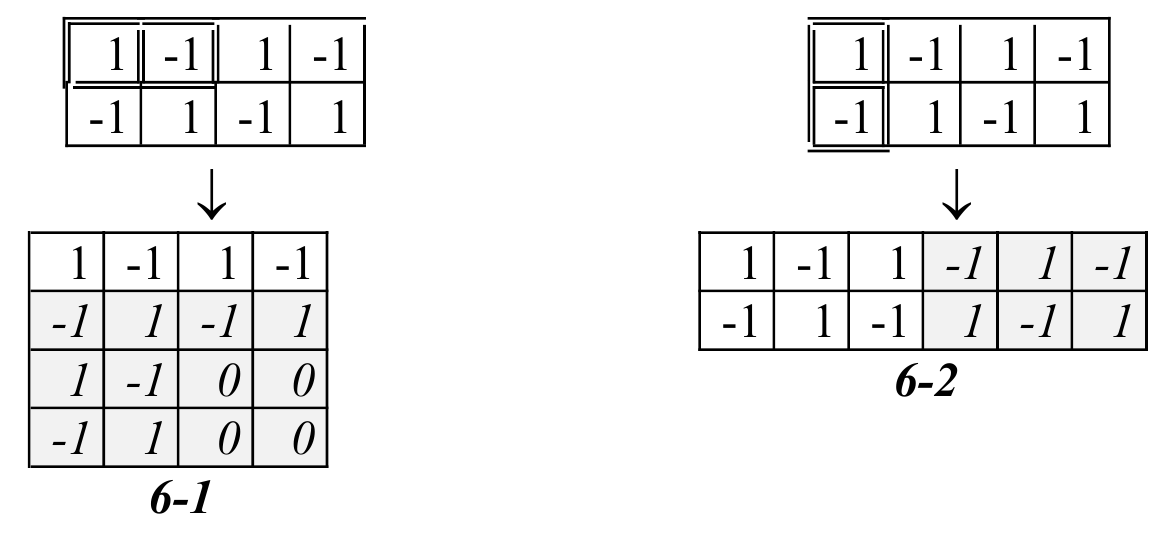}

For each of the cases, the x -distribution is unique.
The corresponding graphs look like this:
\begin{figure}[ht] 
\center{\includegraphics[height=2.0cm]{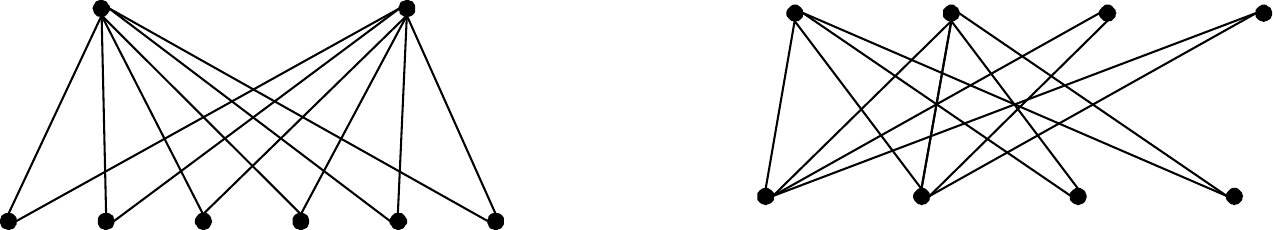}}
\caption{ }
\end{figure}

The corresponding sweeps

\begin{figure}[ht] 
\center{\includegraphics[height=3.5cm]{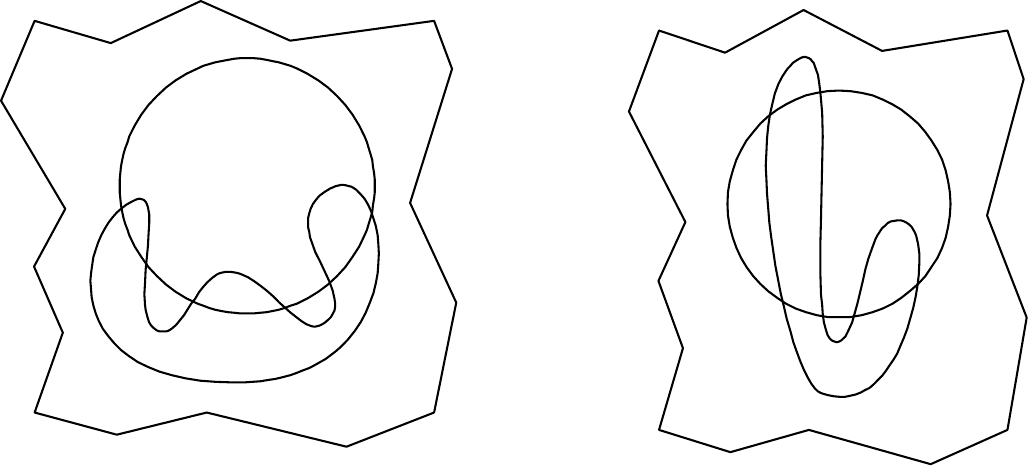}}
\caption{ }
\end{figure}

Thus, the following theorem holds.
\begin{theorem}
There are 2 configurations with 6 intersections: $N( 6 ) = 2$.
\end{theorem} 

 \section{ Eight intersection points}

       \subsection{ Let's consider all intersections of matrix 6-1 }
       
 For all considered for this matrix , there are 4 cases of friction x-distribution is uniform.

 \subsubsection{Case 1}                

\includegraphics[height=2.5cm]{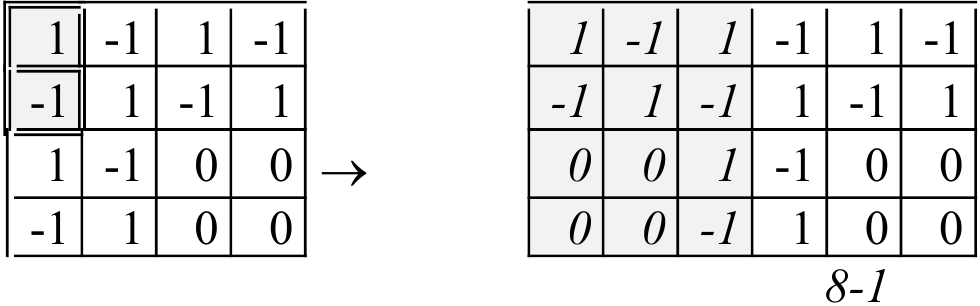}

Corresponding graph and sweep:

\begin{figure}[ht] 
\center{\includegraphics[height=4.5cm]{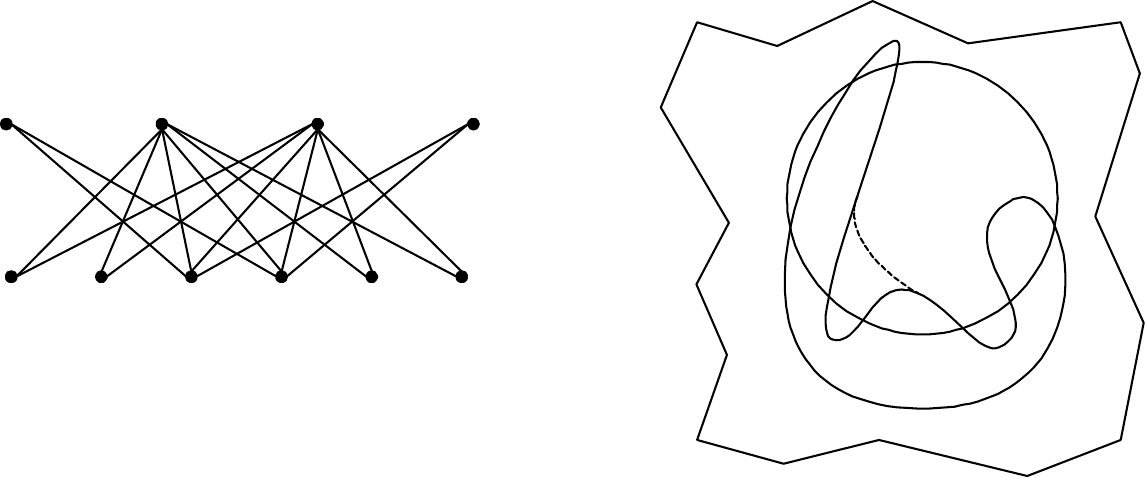}}
\caption{ }
\end{figure}
The dotted line shows the sweep before stretching.

\subsubsection{Case 2}
      \medskip
       
\includegraphics[height=2.5cm]{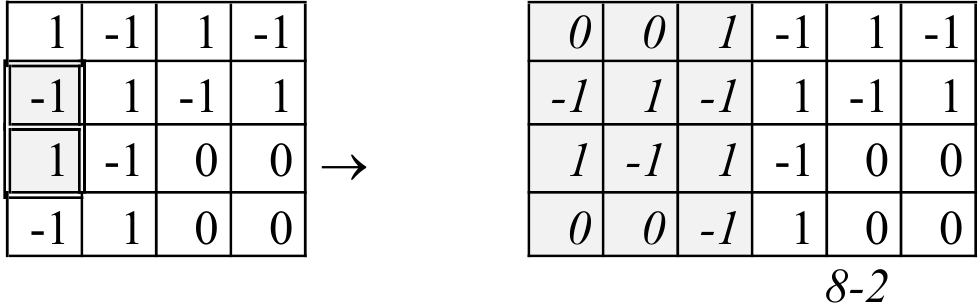}

Corresponding graphs and examples of sweeps of one configuration:

\begin{figure}[ht] 
\center{\includegraphics[height=4.5cm]{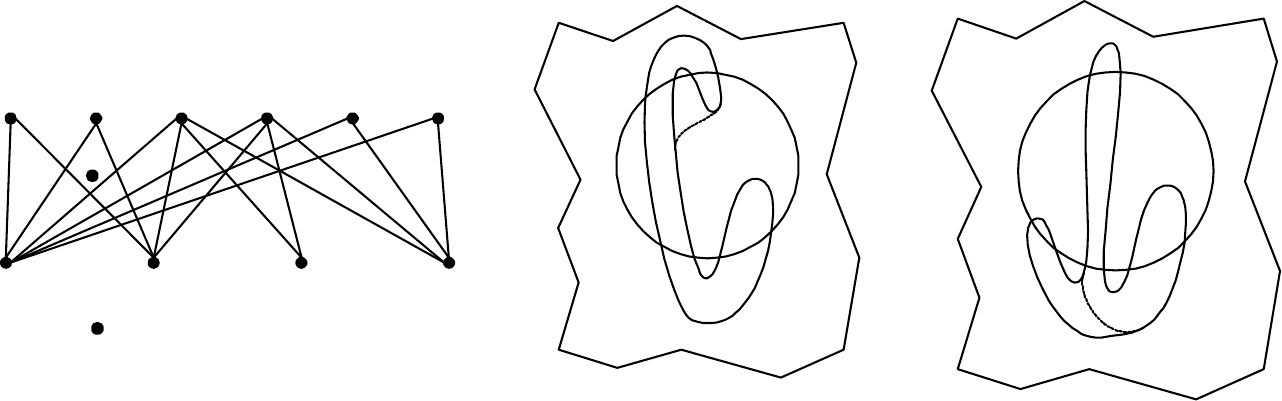}}
\caption{ }
\end{figure}

\newpage

 \subsubsection{ Case 3}

     \medskip
       
\includegraphics[height=2.5cm]{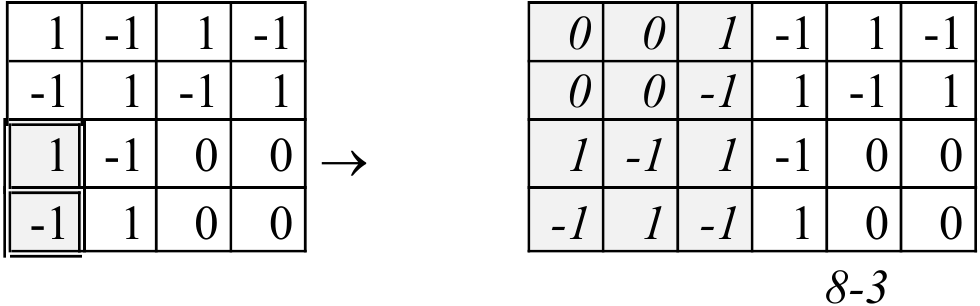}

Corresponding graphs and examples of sweeps of one configuration:

\begin{figure}[ht] 
\center{\includegraphics[height=4.5cm]{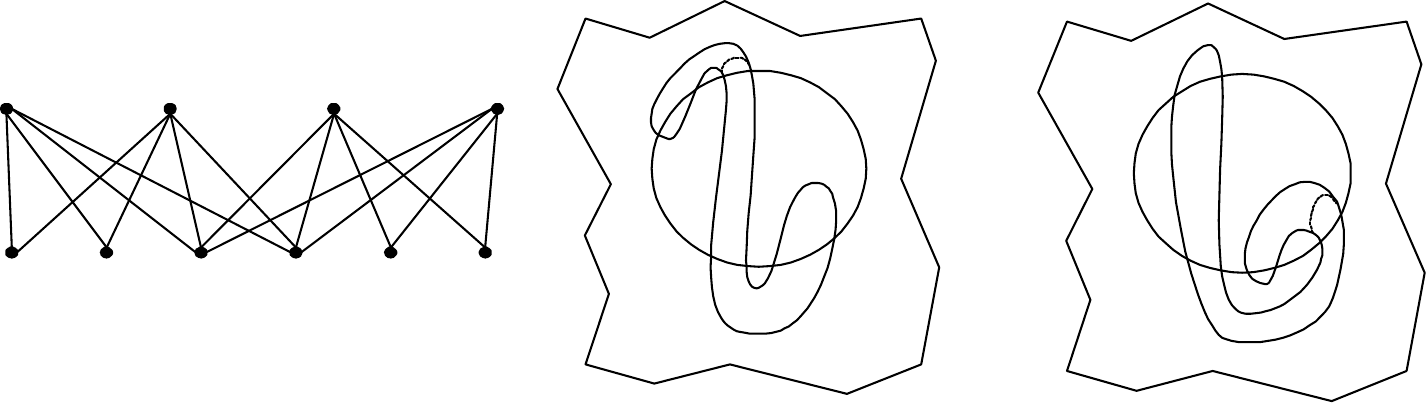}}
\caption{ }
\end{figure}

\subsubsection{Case 4}

  \medskip
       
\includegraphics[height=2.5cm]{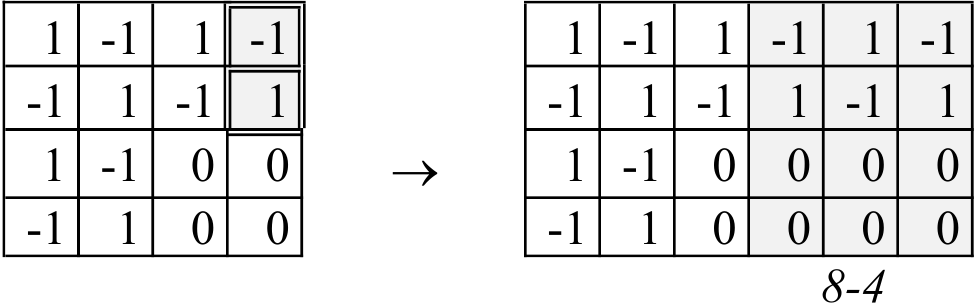}

Corresponding graphs and examples of sweeps of one configuration:

Matrix 8-4 is equal to matrix 8-1 up to the permutation of the columns, so this case does not add new configurations.

\subsection{Now consider the cross-section of matrix 6-2.}              
\subsubsection{ Case 1}

 \medskip
       
\includegraphics[height=1.5cm]{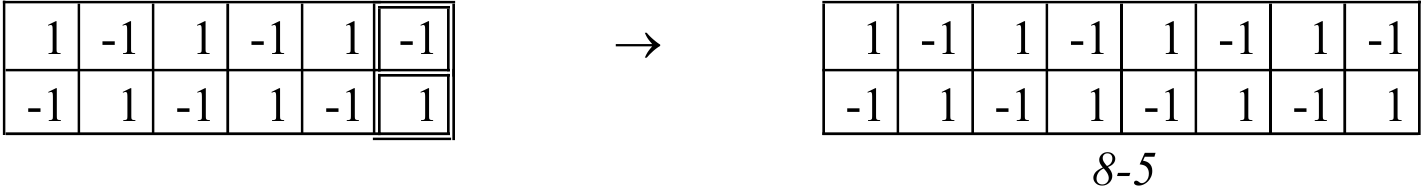}

Corresponding graphs and examples of sweeps of one configuration:

\begin{figure}[ht] 
\center{\includegraphics[height=4.5cm]{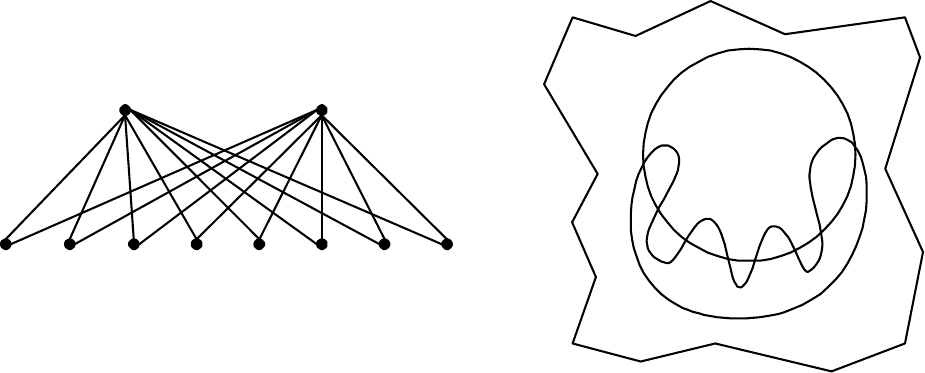}}
\caption{ }
\end{figure}

It is clear that in this case the x-distribution is unique.

 \subsubsection{Case 2}
 
For this case, the x-distribution is not unique.
                      
                       a)
\includegraphics[height=2.5cm]{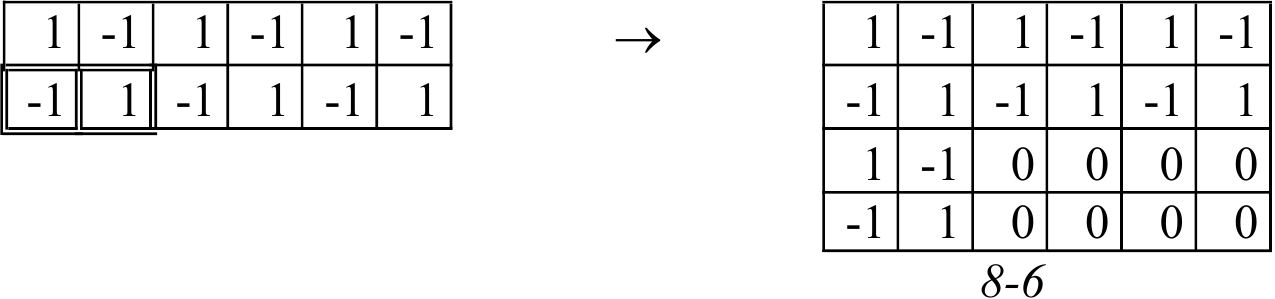}

Corresponding graphs and examples of sweeps of one configuration:

\begin{figure}[ht] 
\center{\includegraphics[height=4.5cm]{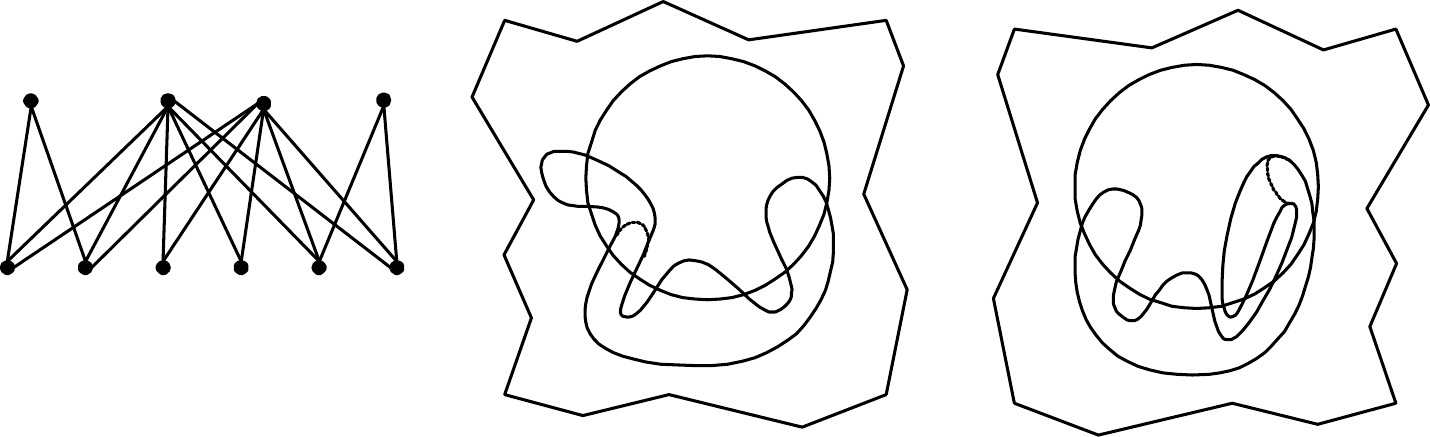}}
\caption{ }
\end{figure}

Corresponding graphs and examples of sweeps of one configuration:

                       b)

\includegraphics[height=2.5cm]{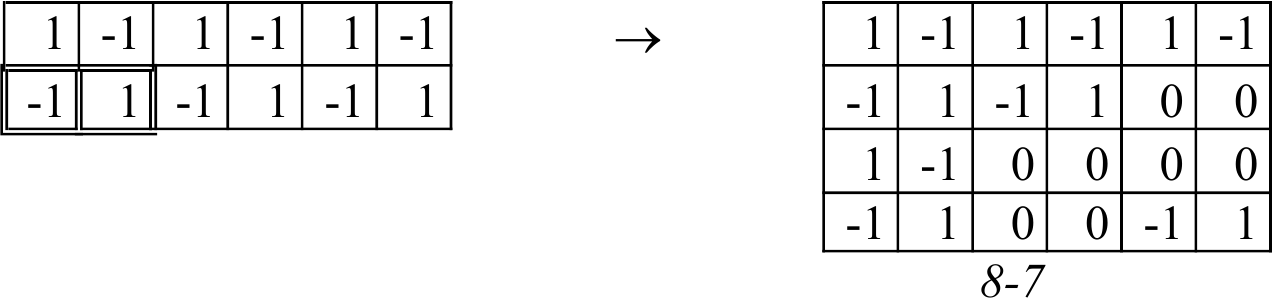}

Corresponding graphs and examples of sweeps of one configuration:
\newpage

\begin{figure}[ht] 
\center{\includegraphics[height=4.5cm]{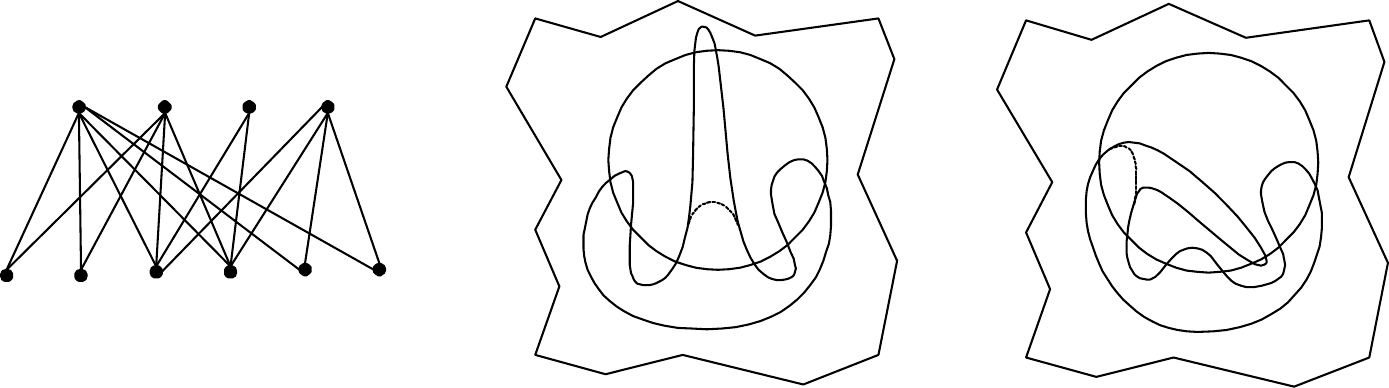}}
\caption{ }
\end{figure}

Matrices 8-6 and 8-7 do not form new configurations, because they are equal to matrices 8-1 and 8-2 up to the permutation of the columns.
              
 Now we can write out all possible configurations for 8 intersection points:

\includegraphics[height=2.8cm]{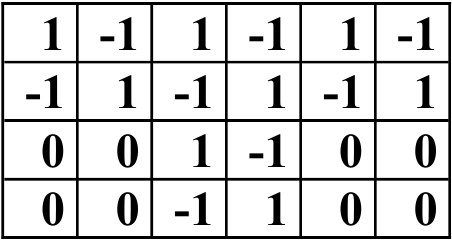}

matrix 8-1
\medskip

\includegraphics[height=7.0cm]{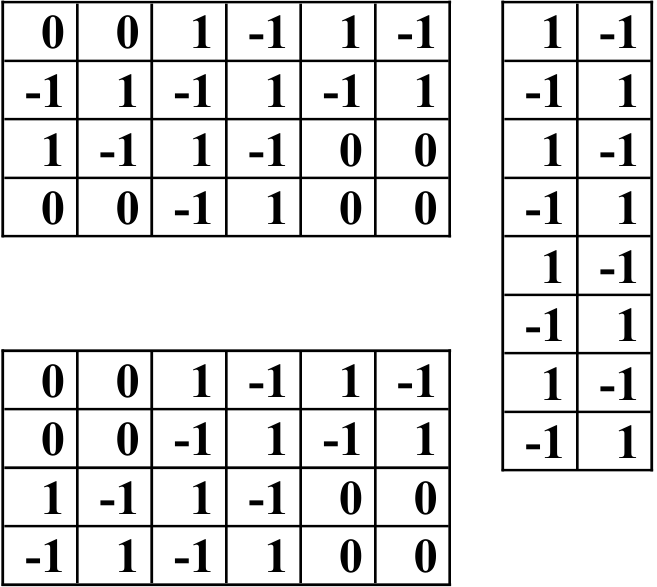}

matrix 8-2, matrix 8-3 and matrix 8-5 (reversed)

\newpage

\begin{figure}[ht] 
\center{\includegraphics[height=4.1cm]{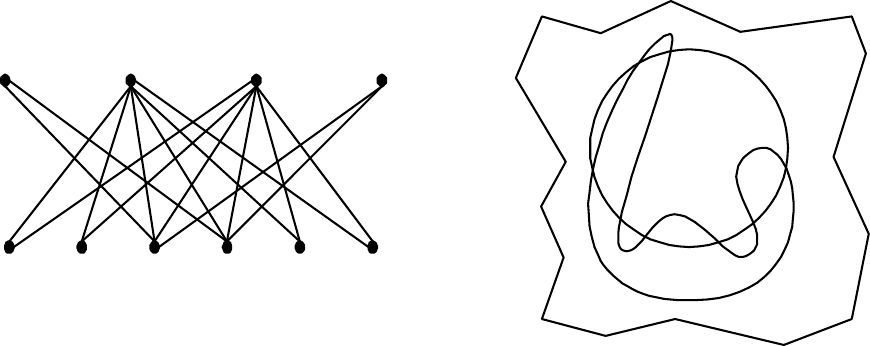},
\includegraphics[height=4.1cm]{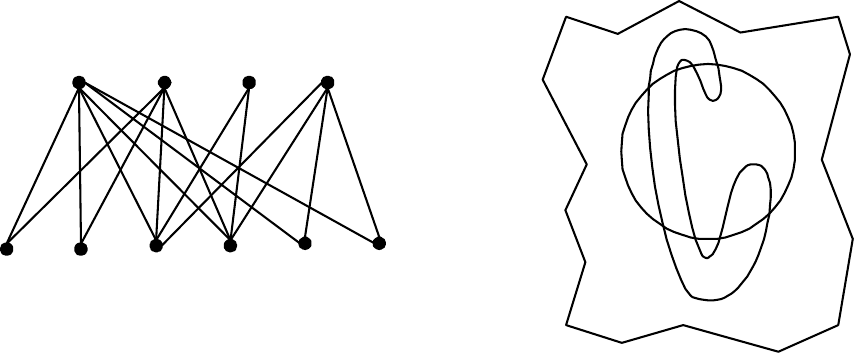},
\includegraphics[height=4.1cm]{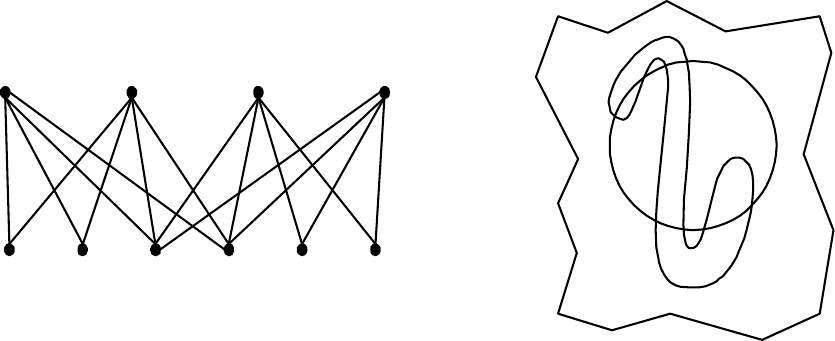},
\includegraphics[height=4.1cm]{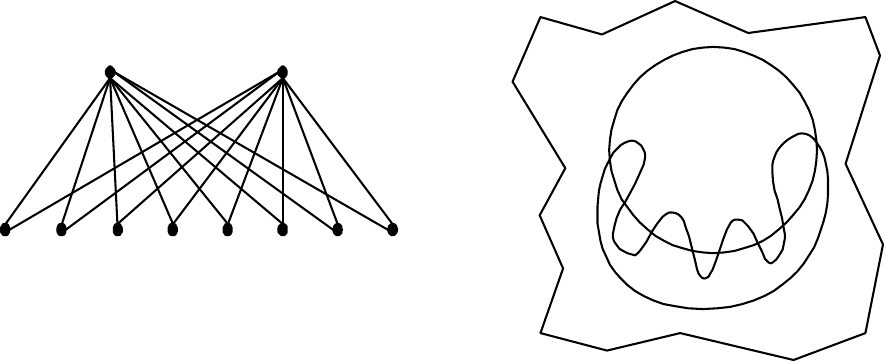}}
\caption{ All possible configurations with 8 intersections }
\end{figure}

Thus, the following theorem holds.
\begin{theorem}
There are 4 configurations with 8 intersections: $N( 8 ) = 4$.
\end{theorem}

\newpage

\section{ Ten intersections}

Later, we  use the determinant vectors of the matrices in order to distinguish or check the equality of the studied matrices.

    \subsection{We consider the intersection of matrix 8-1}
    \subsubsection{Case 1}             
x -distribution is uniform.

\includegraphics[height=6.5cm]{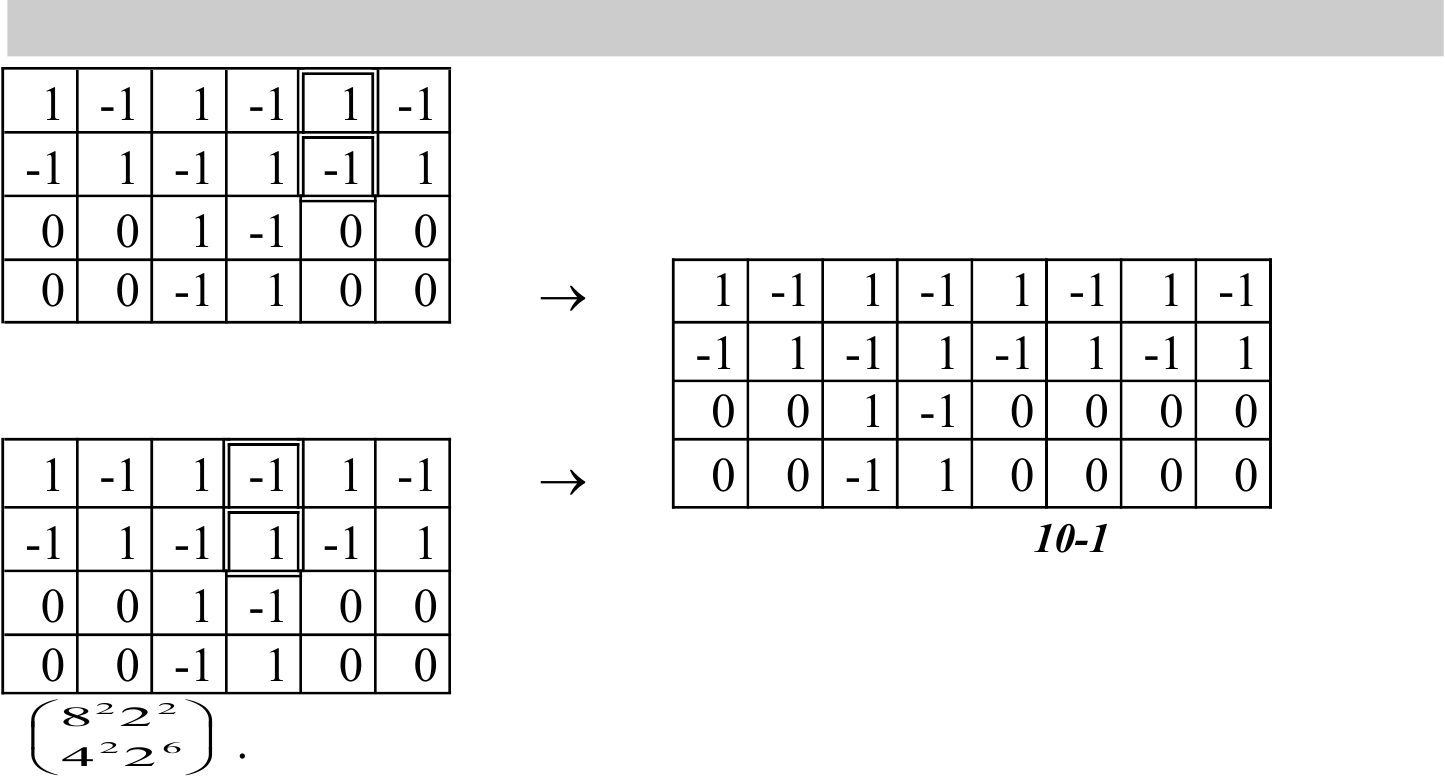}

The initial vectors and its  matrix.

\begin{figure}[ht] 
\center{\includegraphics[height=3.5cm]{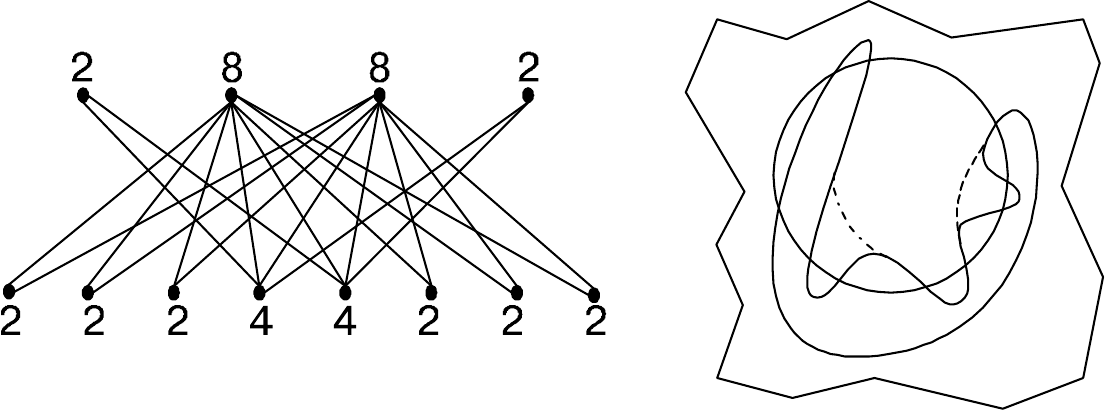}}
\label{p1}
\caption{ }
\end{figure}

\subsubsection{Case 2}             
x -distribution is uniform.

\includegraphics[height=3.5cm]{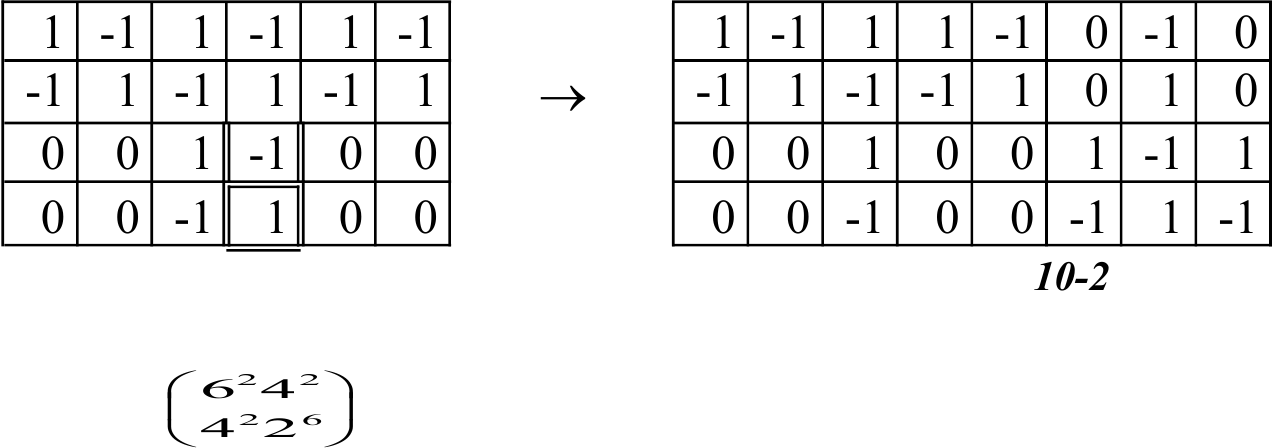}

The initial vectors and its  matrix. Corresponding Count and sweep:
\newpage

\begin{figure}[ht] 
\center{\includegraphics[height=3.5cm]{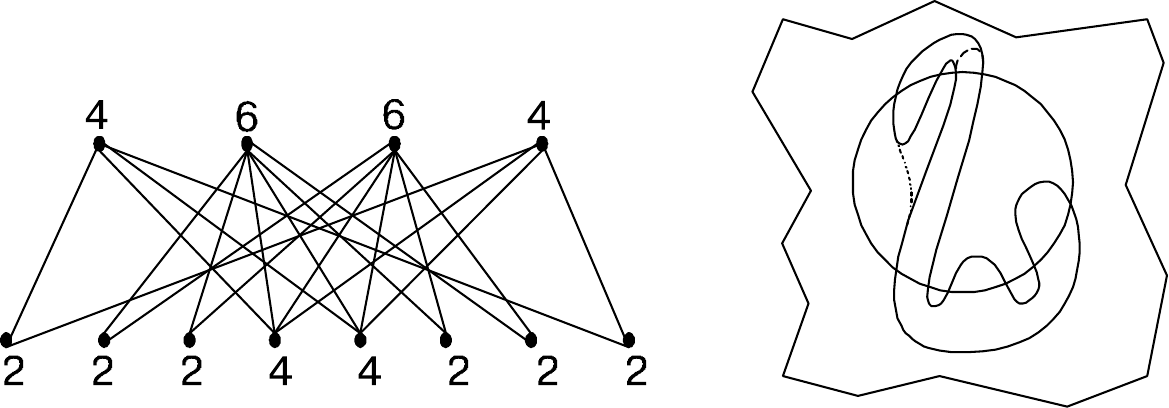}}
\label{p2}
\caption{ }
\end{figure}

 \subsubsection{Case 3}             
x -distribution is uniform.

\includegraphics[height=3.0cm]{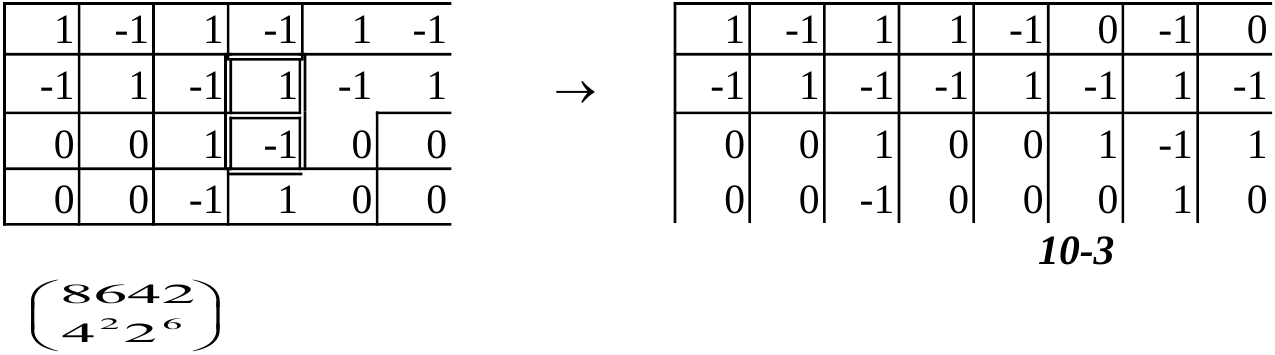}

The initial vectors and its  matrix. Corresponding Count and sweep:

\begin{figure}[ht] 
\center{\includegraphics[height=3.5cm]{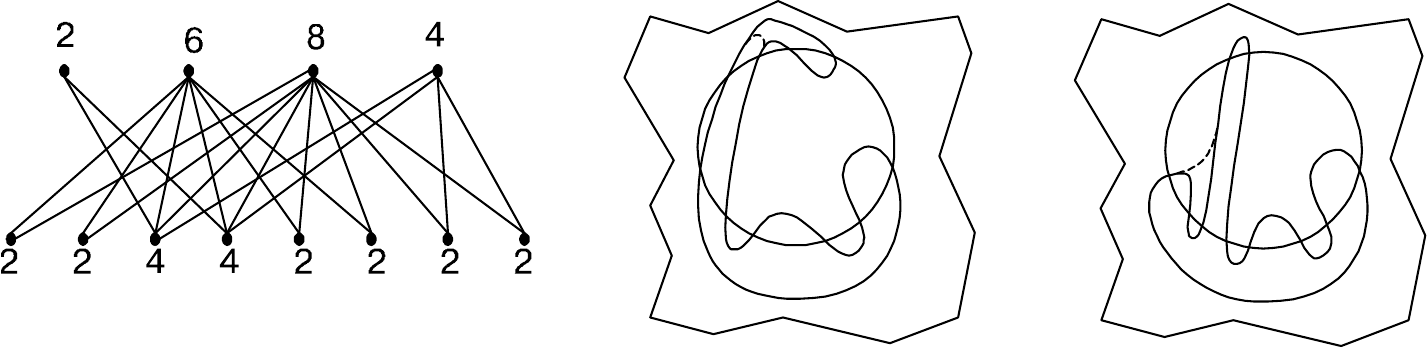}}
\label{p3}
\caption{ }
\end{figure}

 \subsubsection{Case 4}             
x -distribution is not uniform.

a) 

\includegraphics[height=4.5cm]{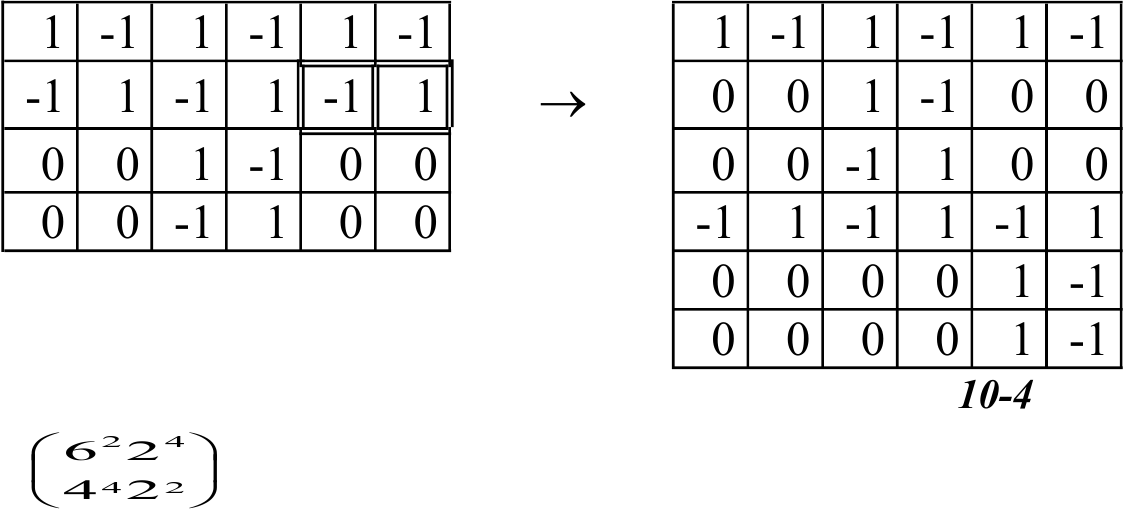}

The initial vectors and its  matrix. Corresponding Count and sample sweeps.

\begin{figure}[ht] 
\center{\includegraphics[height=3.5cm]{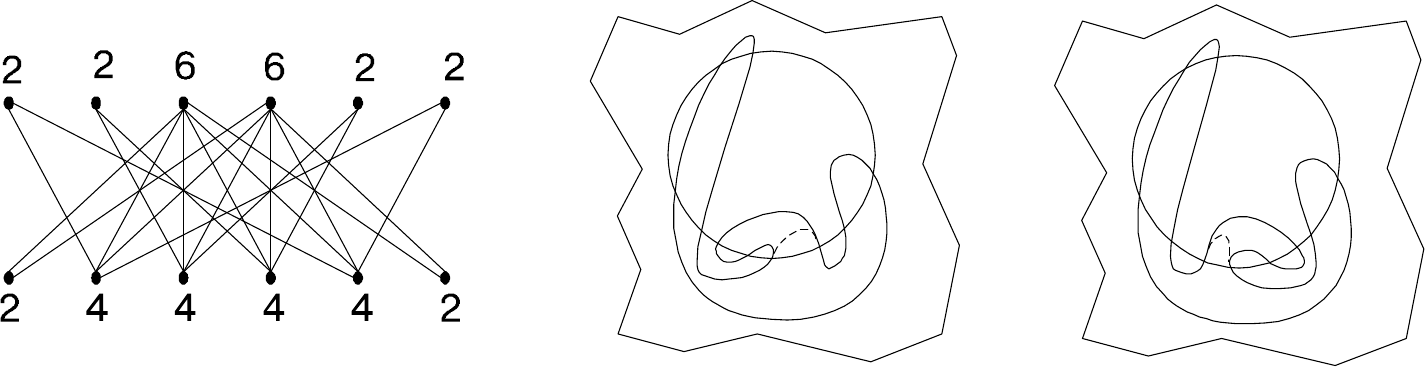}}
\label{p4}
\caption{ }
\end{figure}

\newpage

                       b)

\includegraphics[height=4.5cm]{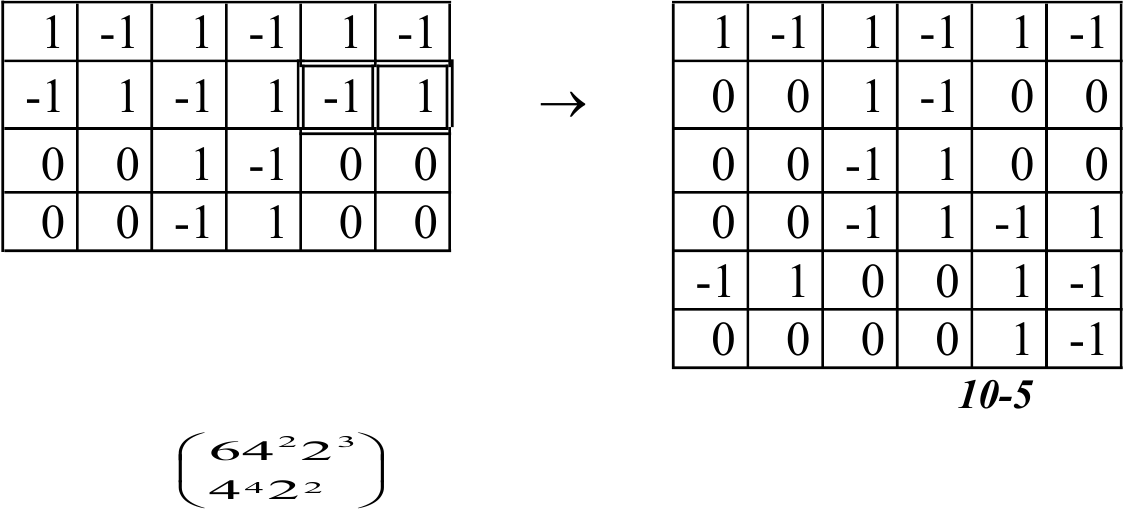}

The initial vectors and its  matrix. Corresponding Count and an example sweep.

\begin{figure}[ht] 
\center{\includegraphics[height=3.5cm]{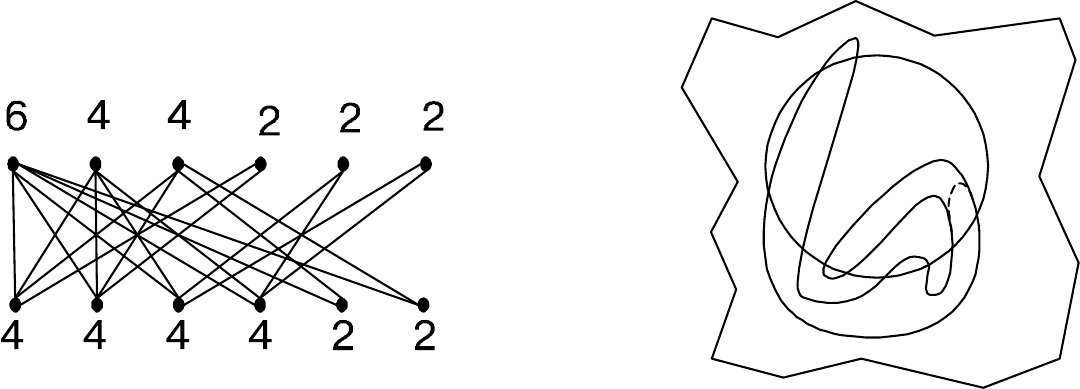}}
\label{p5}
\caption{ }
\end{figure}

\subsubsection{Case 5}

\includegraphics[height=4.5cm]{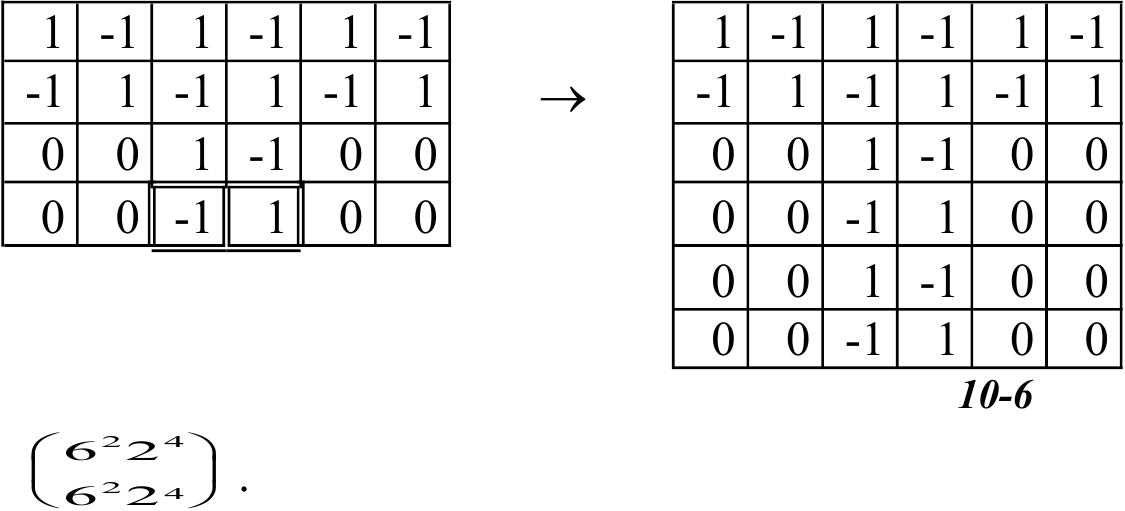}

The initial vectors and its  matrix. Corresponding Count and an example sweep.

\begin{figure}[ht] 
\center{\includegraphics[height=3.5cm]{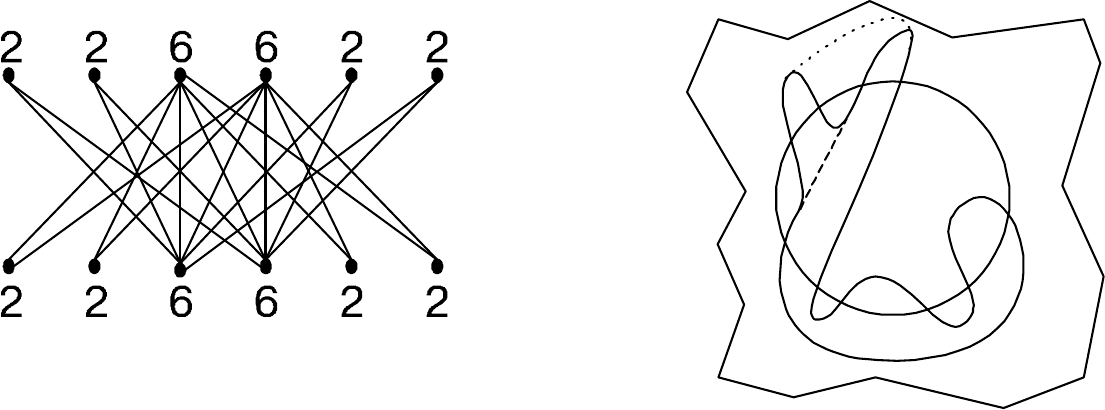}}
\label{p6}
\caption{ }
\end{figure}

 \subsubsection{Case 6}
x -distribution is not uniform.

                       a )
 \medskip                      
                       
\includegraphics[height=4.5cm]{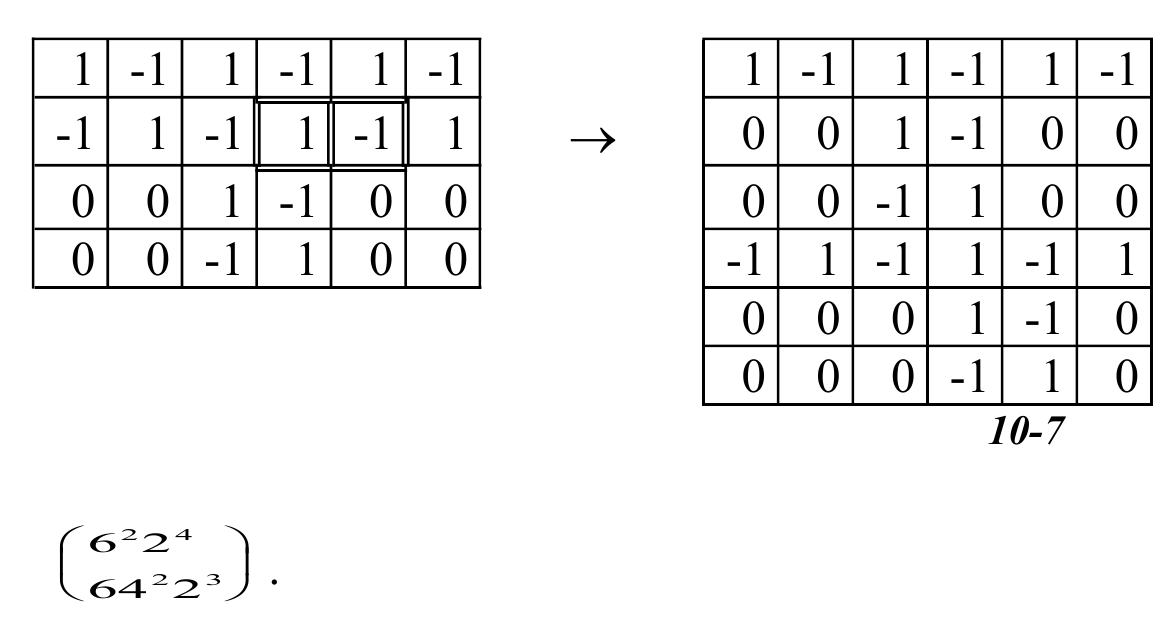}

The initial vectors and its  matrix. Corresponding Count and an example sweep.

\begin{figure}[ht] 
\center{\includegraphics[height=3.5cm]{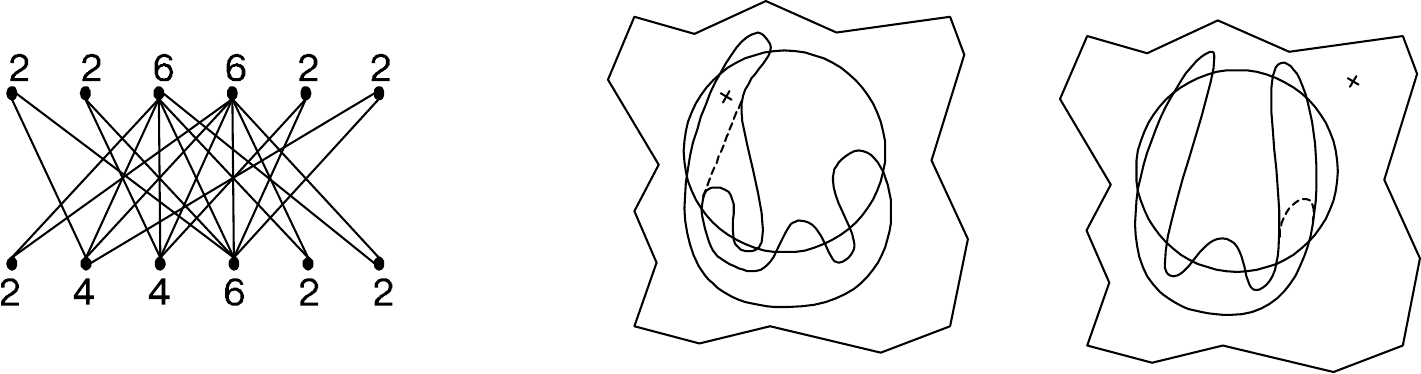}}
\label{p7}
\caption{ }
\end{figure}

                       b )
 \medskip                      
                       
\includegraphics[height=4.5cm]{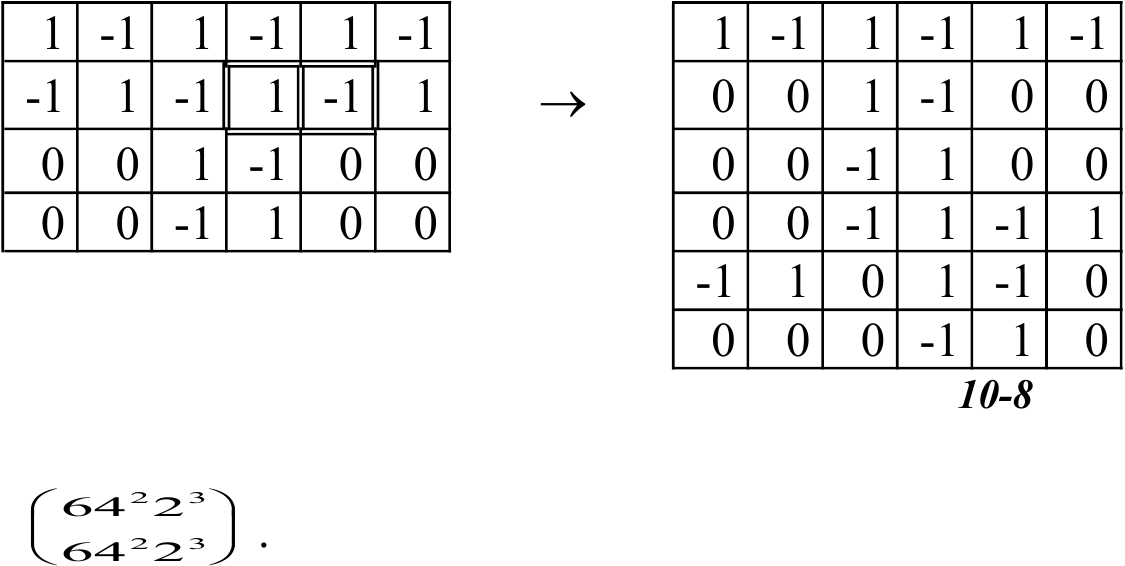}

The initial vectors and its  matrix. Corresponding Count and sample sweeps.

\begin{figure}[ht] 
\center{\includegraphics[height=3.5cm]{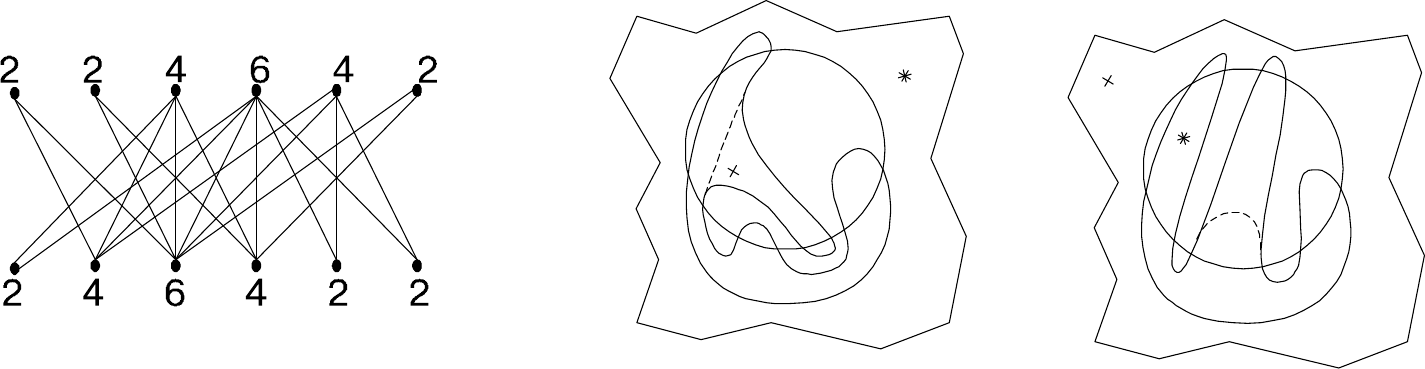}}
\label{p8}
\caption{ }
\end{figure}

            \subsection{ We consider the intersection of matrix 8-3}
                
\subsubsection{ Case 1}
x -distribution is uniform.

\medskip                      
                       
\includegraphics[height=5.5cm]{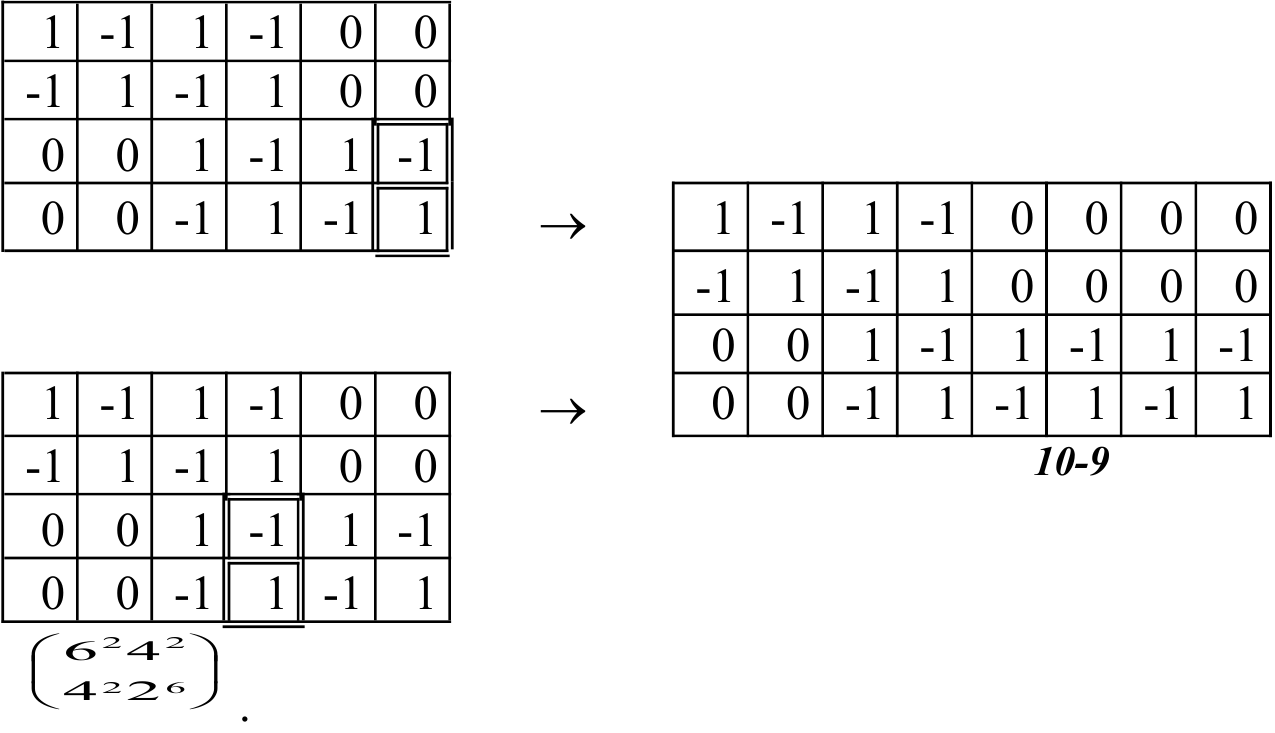}

The initial vectors and its  matrix. Corresponding graph and examples of sweeps.

\begin{figure}[ht] 
\center{\includegraphics[height=3.5cm]{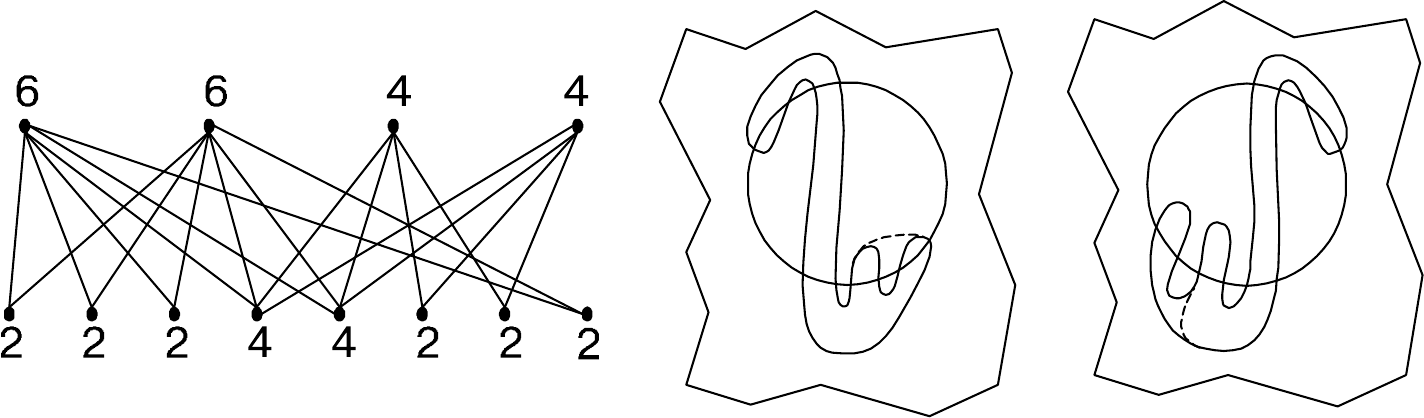}}
\caption{ }
\end{figure}

\subsubsection{ Case 2}
x -distribution is uniform.

\medskip                      
                       
\includegraphics[height=3.5cm]{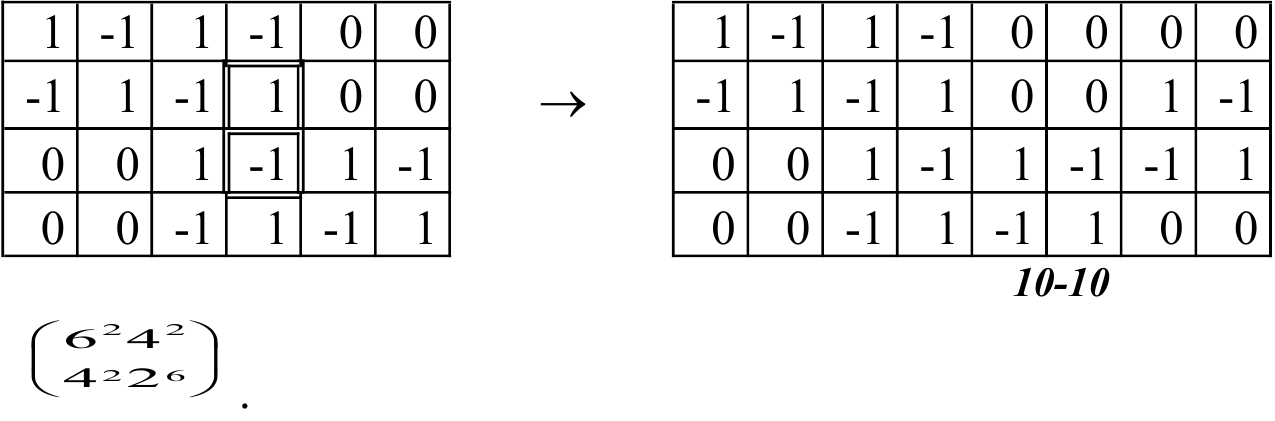}

The initial vectors and its  matrix. Corresponding graph and examples of sweeps.

\begin{figure}[ht] 
\center{\includegraphics[height=3.5cm]{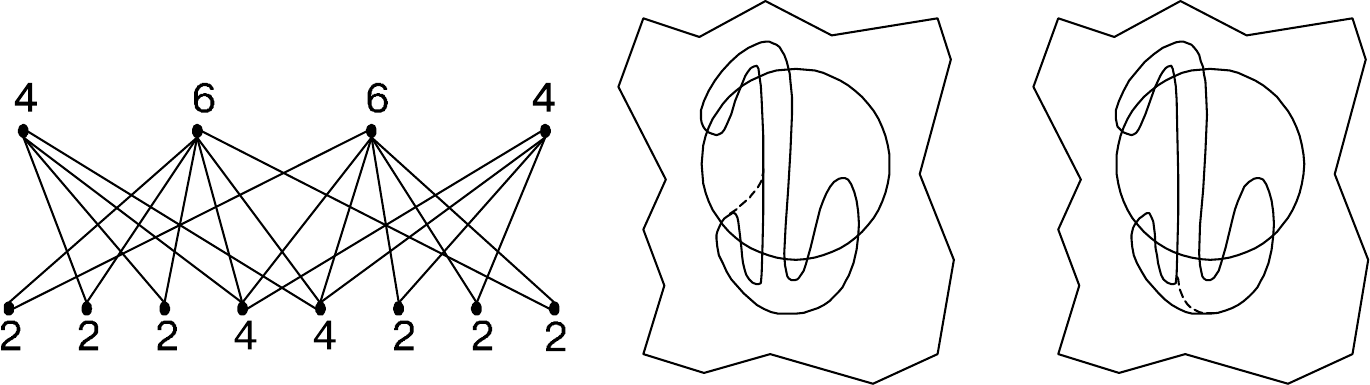}}
\caption{ }
\end{figure}

\subsubsection{ Case 3}
x -distribution is uniform.

\medskip                      
                       
\includegraphics[height=4.5cm]{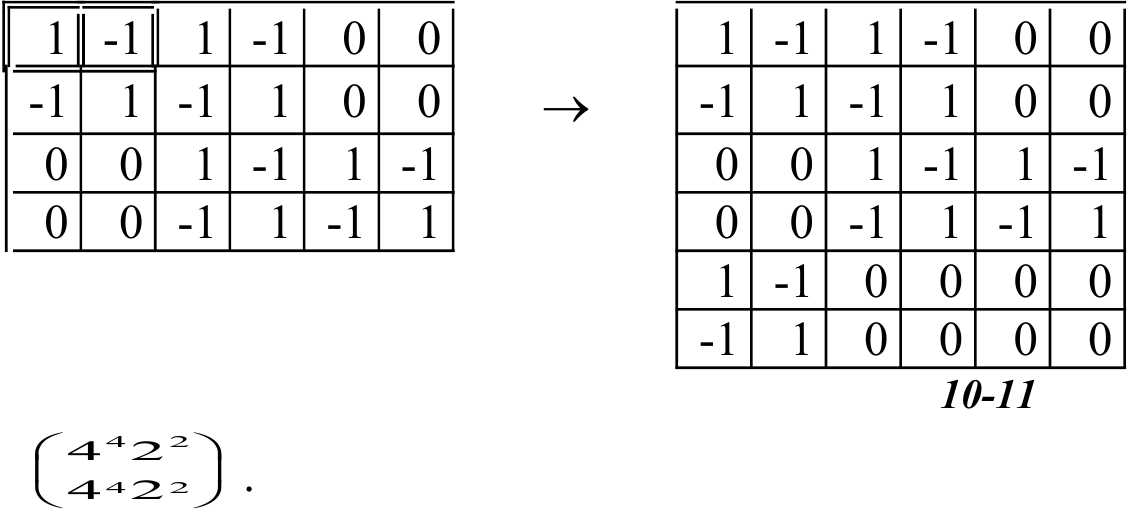}

The initial vectors and its  matrix. Corresponding graph and sample sweeps.

\begin{figure}[ht] 
\center{\includegraphics[height=3.5cm]{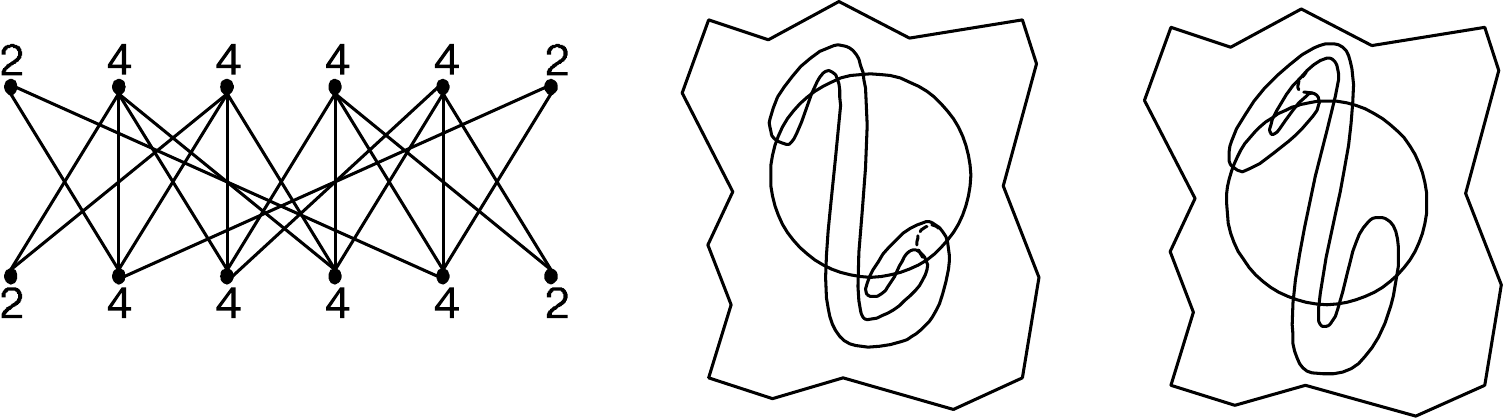}}
\label{p11}
\caption{ }
\end{figure}     

\subsubsection{ Case 4}
x -distribution is uniform.

\medskip                      
                       
\includegraphics[height=4.5cm]{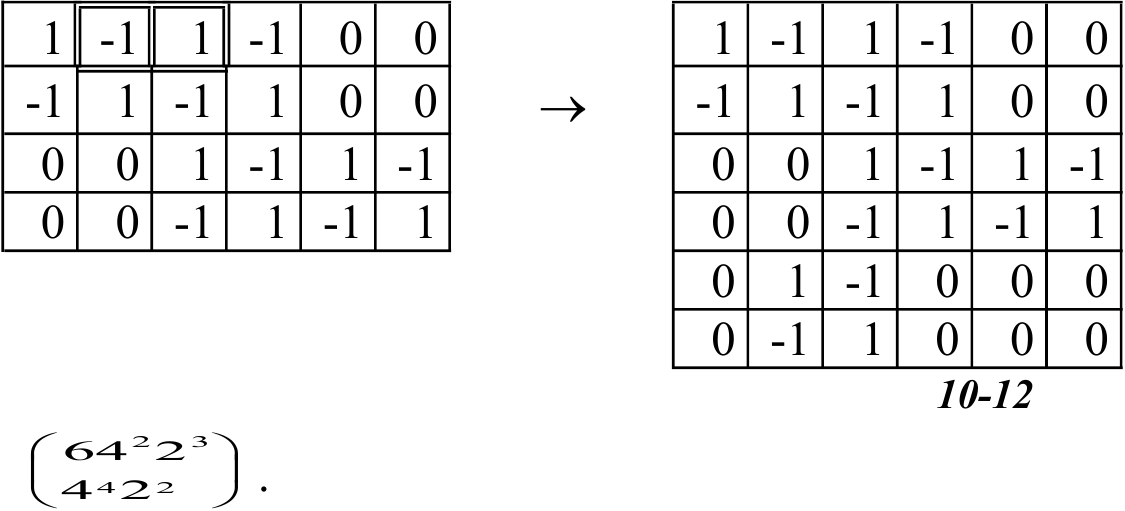}

The initial vectors and its  matrix. Corresponding graph and sample sweeps.

\begin{figure}[ht] 
\center{\includegraphics[height=3.5cm]{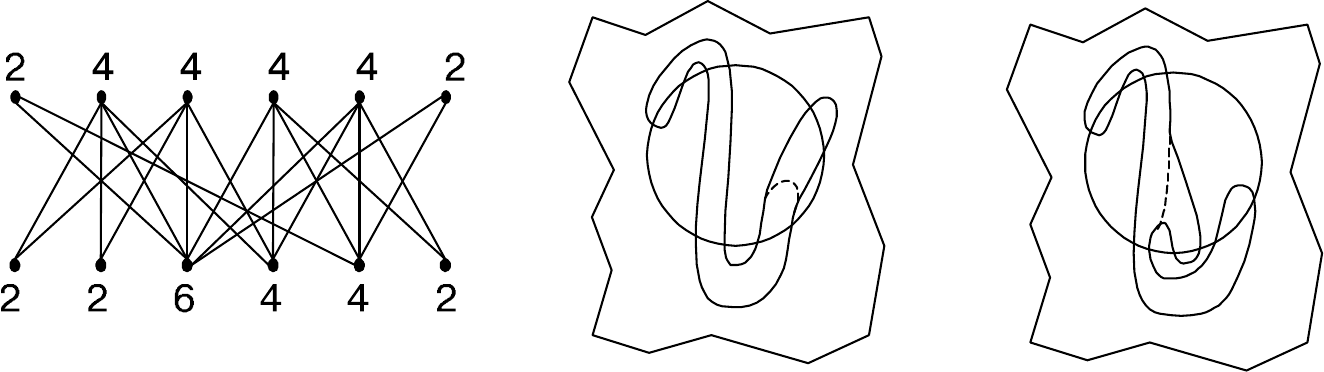}}
\caption{ }
\end{figure}

\subsubsection{ Case 5}
x -distribution is uniform.

\medskip                      
                       
\includegraphics[height=4.5cm]{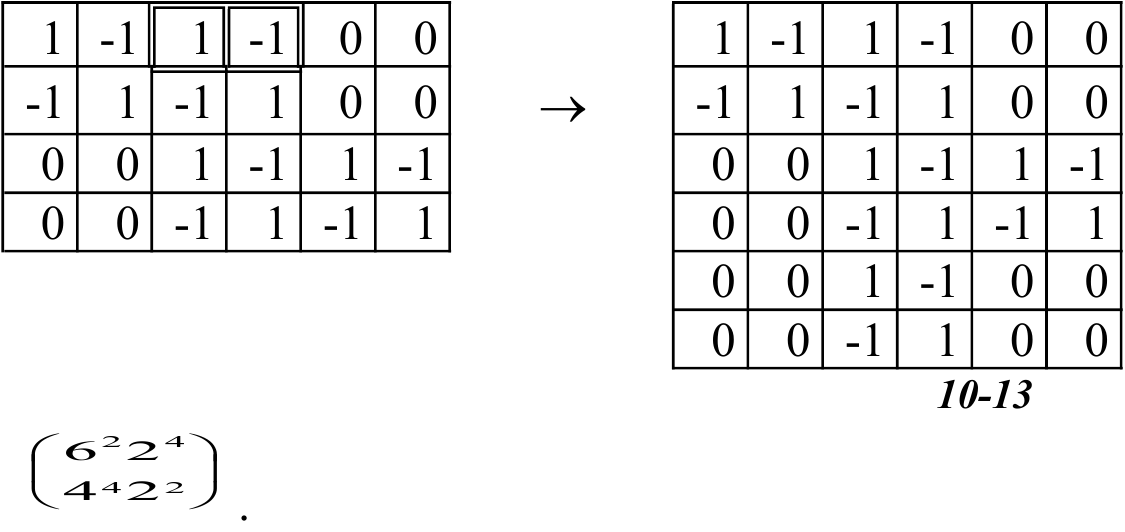}

The initial vectors and its  matrix. Corresponding graph and examples of sweeps.

\begin{figure}[ht] 
\center{\includegraphics[height=3.5cm]{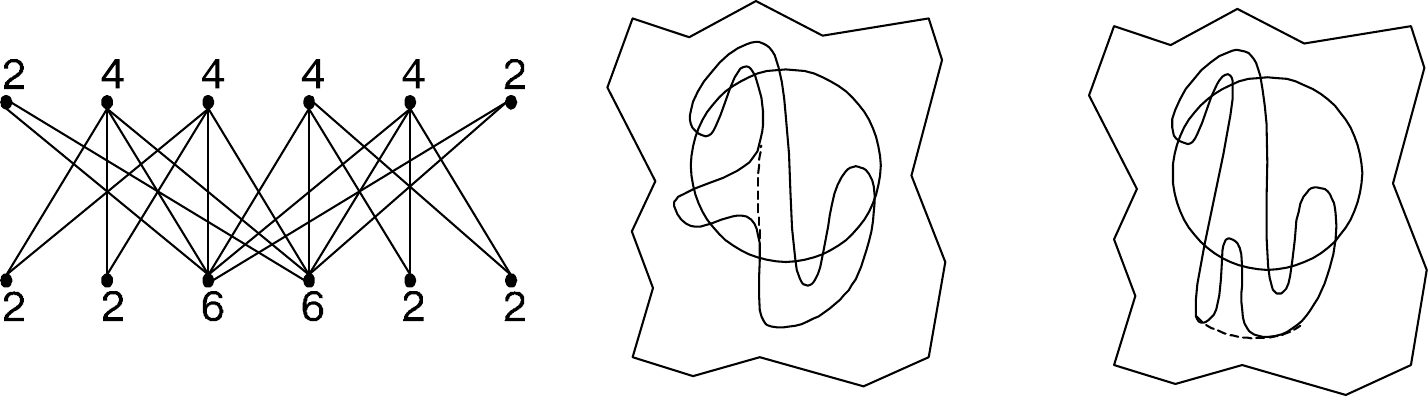}}
\caption{ }
\end{figure}

\subsection{We consider the intersection of the 8-2 matrix}

\subsubsection{ Case 1}
x -distribution is uniform.

\medskip                      
                       
\includegraphics[height=6.5cm]{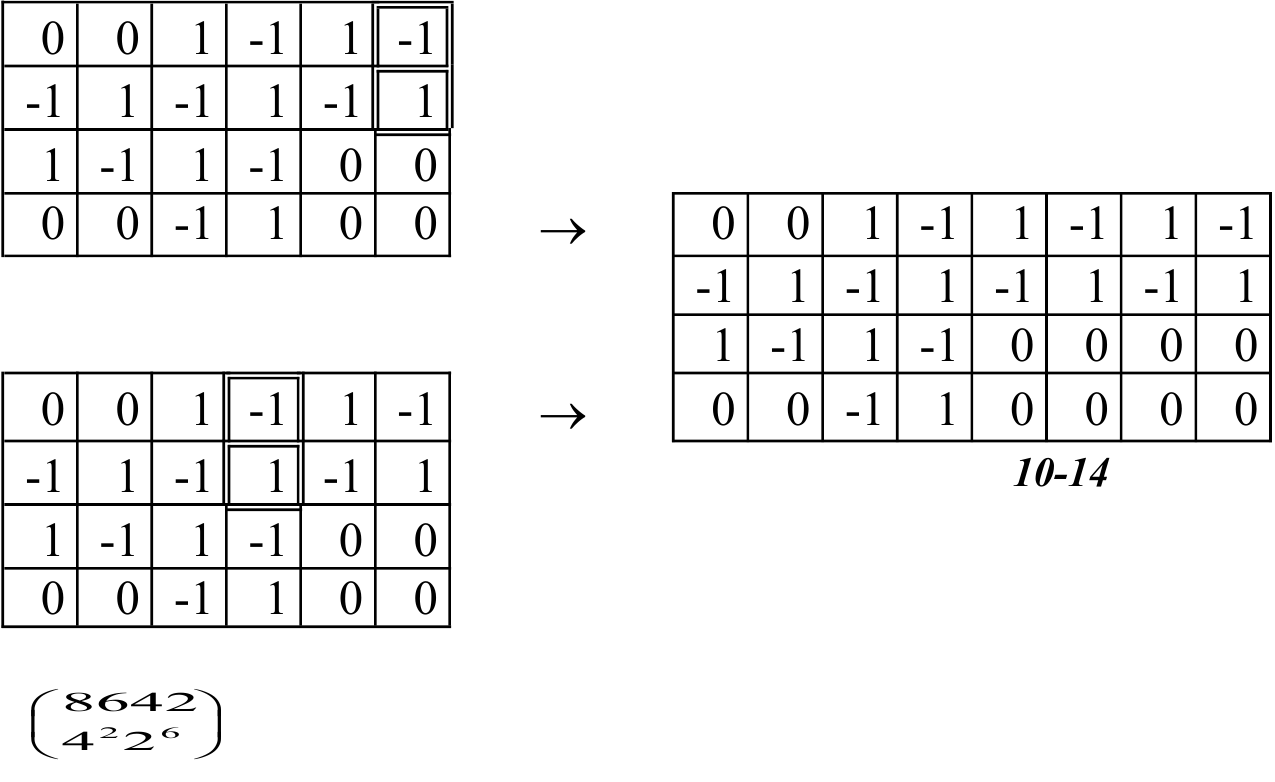}

The initial vectors and its  matrix. Corresponding graph and examples of sweeps.

\begin{figure}[ht] 
\center{\includegraphics[height=3.5cm]{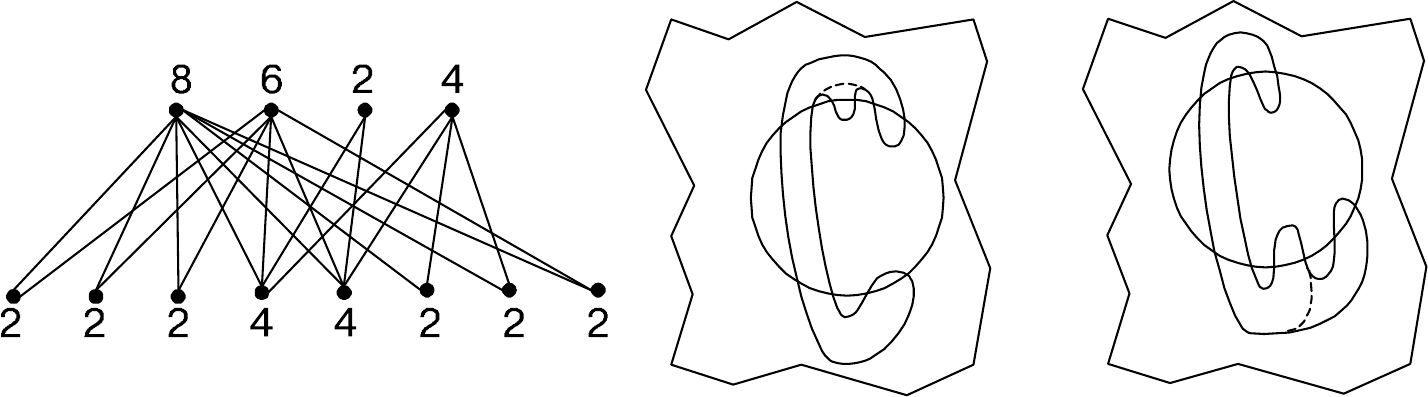}}
\caption{ }
\end{figure}     

 \subsubsection{ Case 2}
x -distribution is uniform.

\medskip                      
                       
\includegraphics[height=3.9cm]{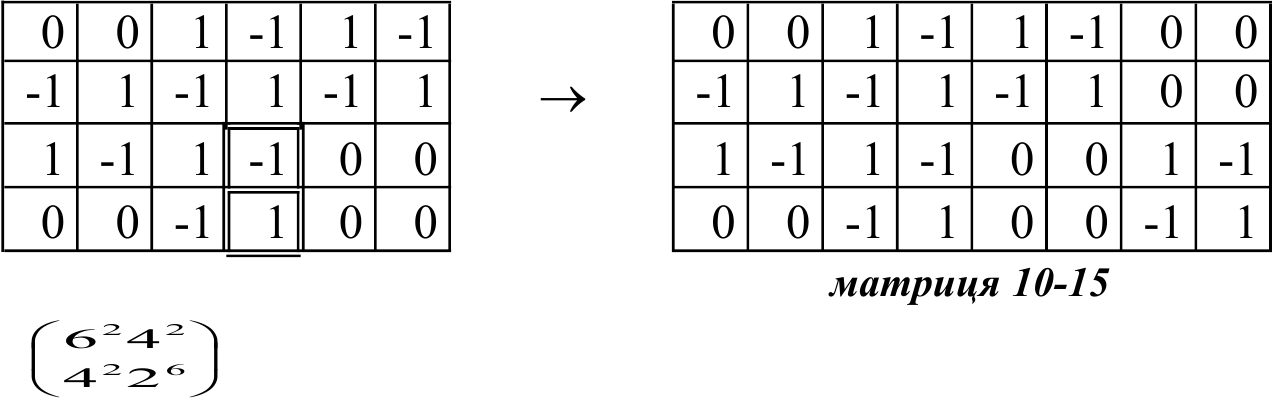}

The initial vectors and its  matrix. Corresponding graph and examples of sweeps.

\begin{figure}[ht] 
\center{\includegraphics[height=3.5cm]{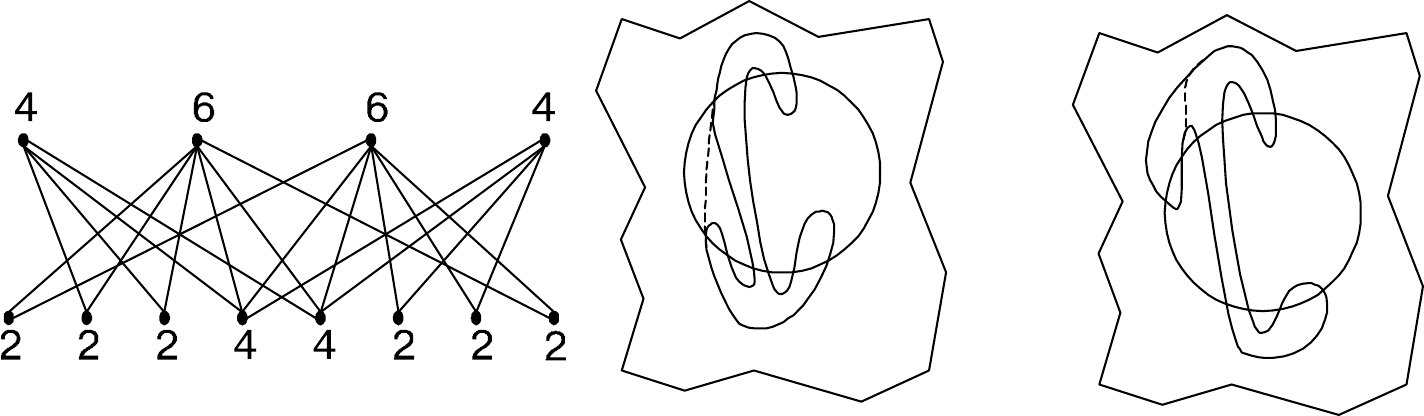}}
\caption{ }
\end{figure} 

 \subsubsection{ Case 3}
x -distribution is uniform.

\medskip                      
                       
\includegraphics[height=6.0cm]{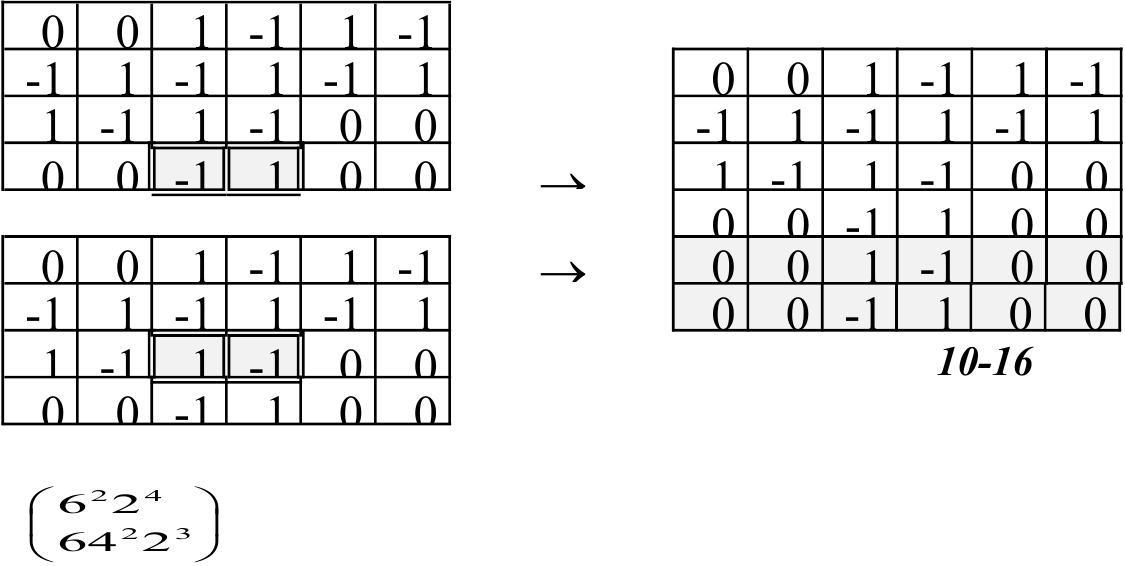}

The initial vectors and its  matrix. Corresponding graph and examples of sweeps.

\begin{figure}[ht] 
\center{\includegraphics[height=3.0cm]{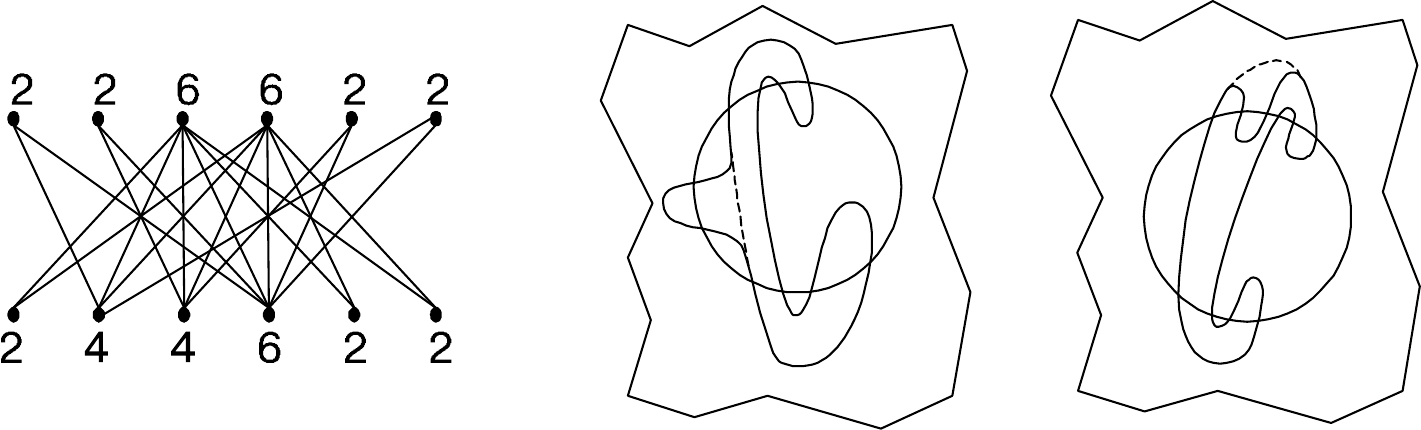}}
\caption{ }
\end{figure} 

               \subsubsection{ Case 4}
x -distribution is uniform.

\medskip                      
                       
\includegraphics[height=4.9cm]{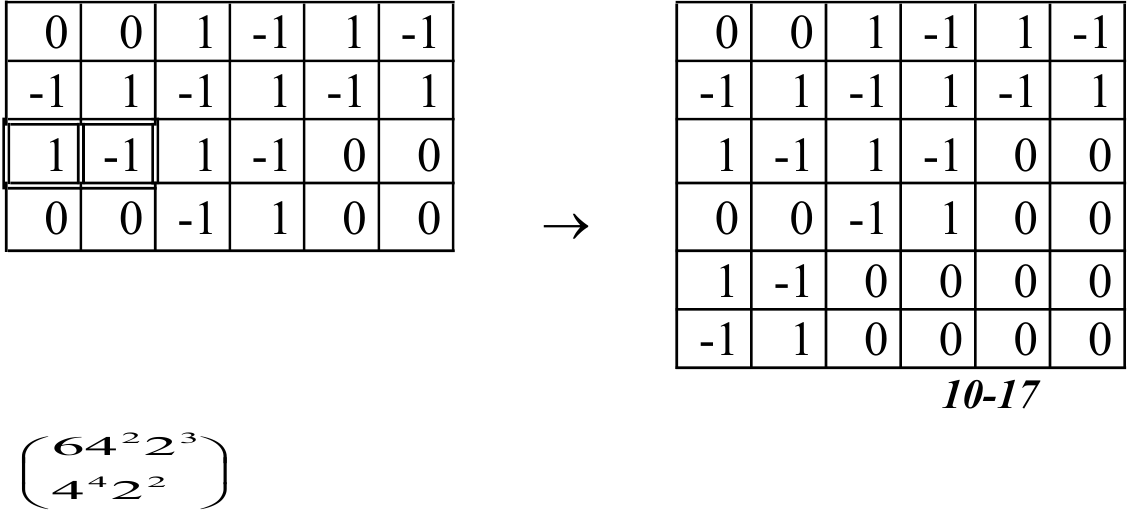}

The initial vectors and its  matrix. Corresponding graph and examples of sweeps.

\begin{figure}[ht] 
\center{\includegraphics[height=3.5cm]{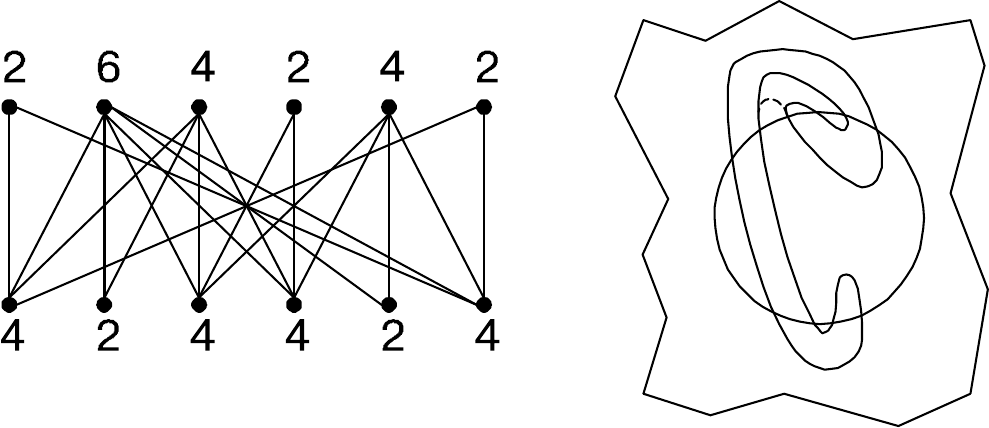}}
\caption{ }
\end{figure} 

               \subsubsection{ Case 5}
x -distribution is uniform.

\medskip                      
                       
\includegraphics[height=3.5cm]{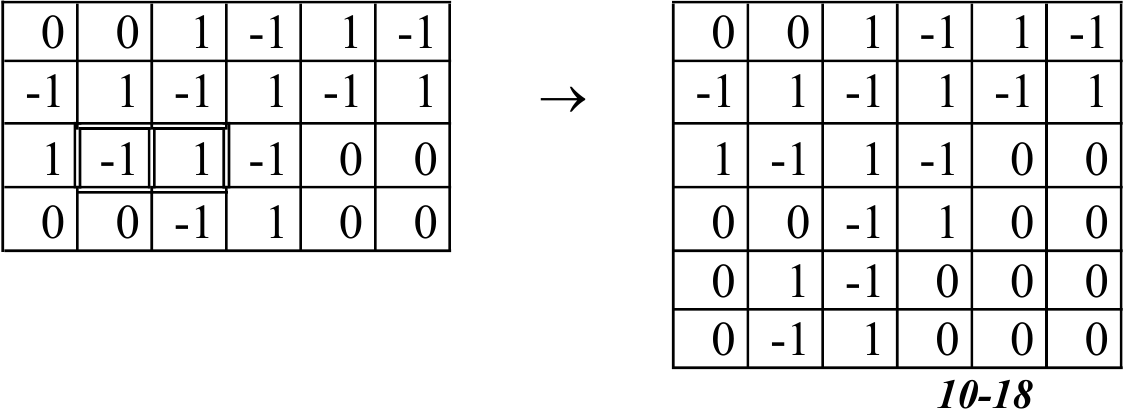}

The initial vectors. Corresponding graph and examples of sweeps.

\begin{figure}[ht] 
\center{\includegraphics[height=3.5cm]{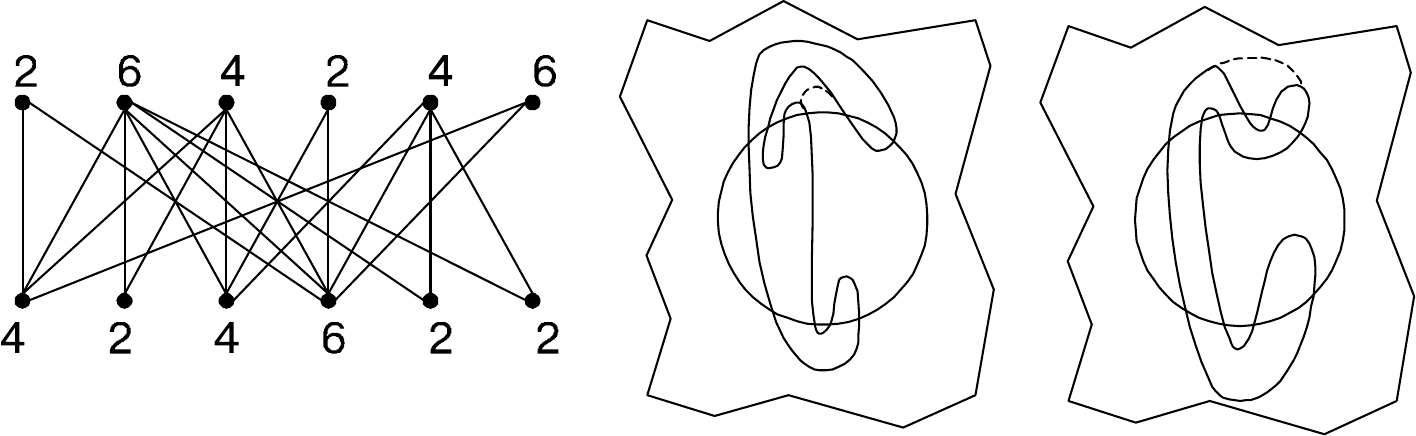}}
\label{p18}
\caption{ }
\end{figure}

            \subsubsection{ Case 6}
x -distribution is not uniform.

                       a)
\medskip                      
                       
\includegraphics[height=3.5cm]{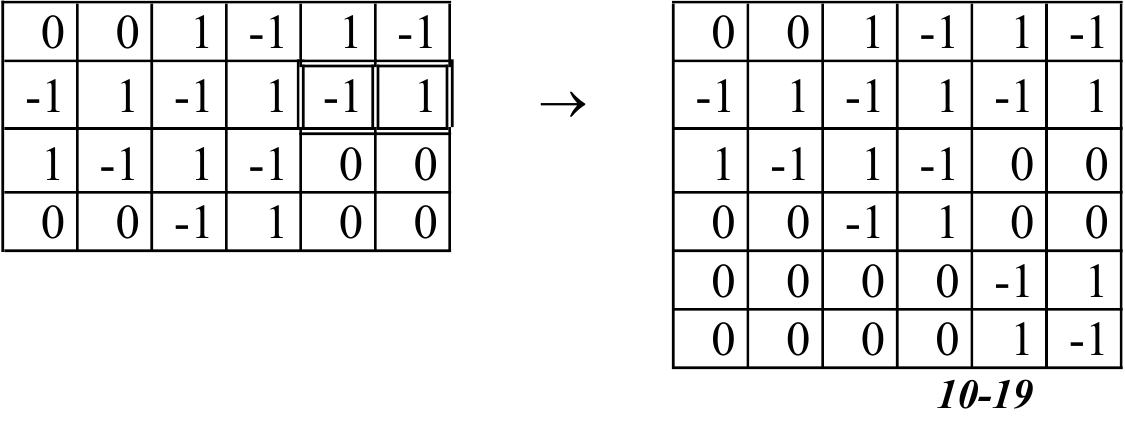}

The initial vectors. Corresponding graph and examples of sweeps.

\begin{figure}[ht] 
\center{\includegraphics[height=3.5cm]{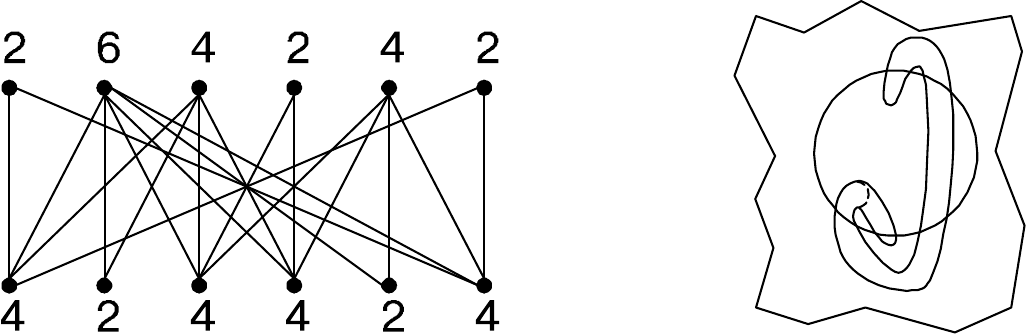}}
\caption{ }
\end{figure}

                       b)
\medskip                      
                       
\includegraphics[height=3.5cm]{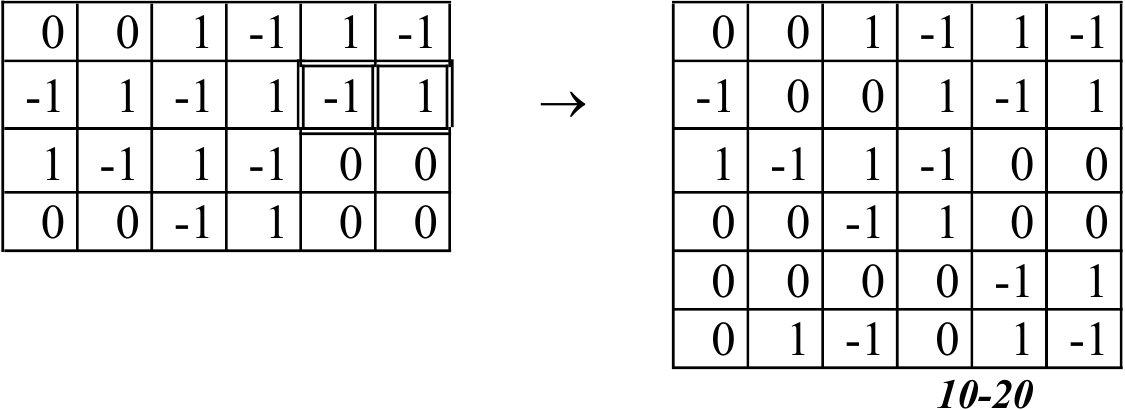}

The initial vectors. Corresponding graph and examples of sweeps.

\begin{figure}[ht] 
\center{\includegraphics[height=3.5cm]{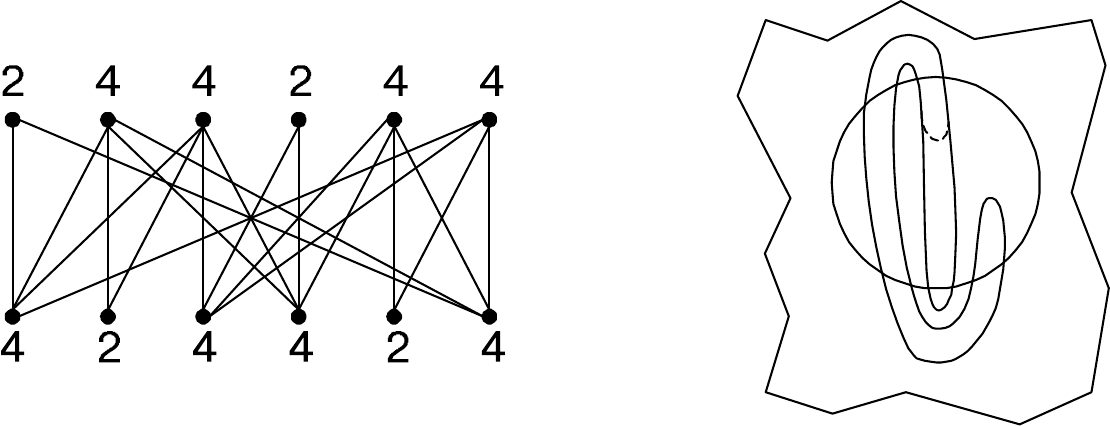}}
\label{p20}
\caption{ }
\end{figure}

               \subsubsection{ Case 7}
x -distribution is uniform.

\medskip                      
                       
\includegraphics[height=3.5cm]{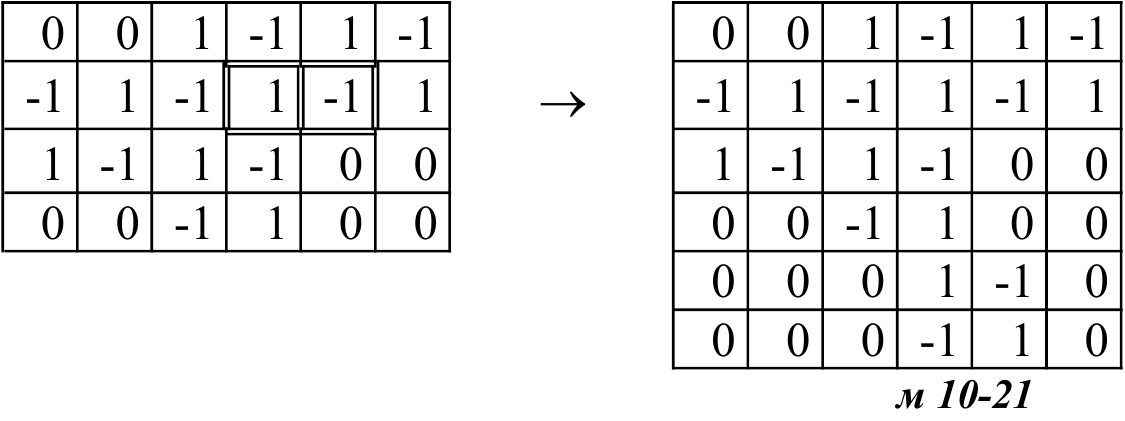}

The initial vectors. Corresponding graph and examples of sweeps.

\begin{figure}[ht] 
\center{\includegraphics[height=3.5cm]{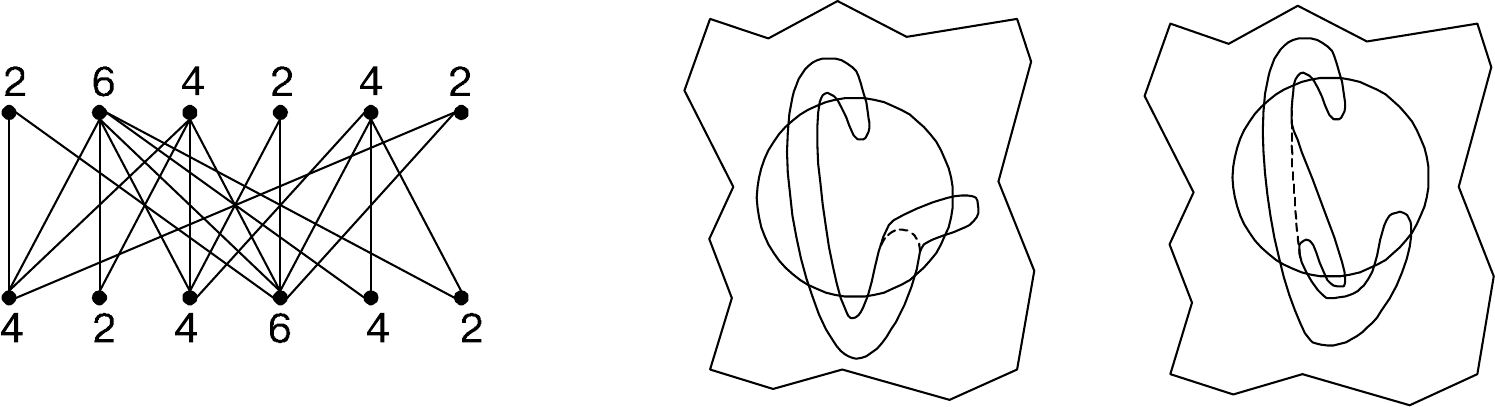}}
\caption{ }
\end{figure} 

\newpage

    \subsubsection{ Case 8}
x -distribution is uniform.

\medskip                      
                       
\includegraphics[height=3.5cm]{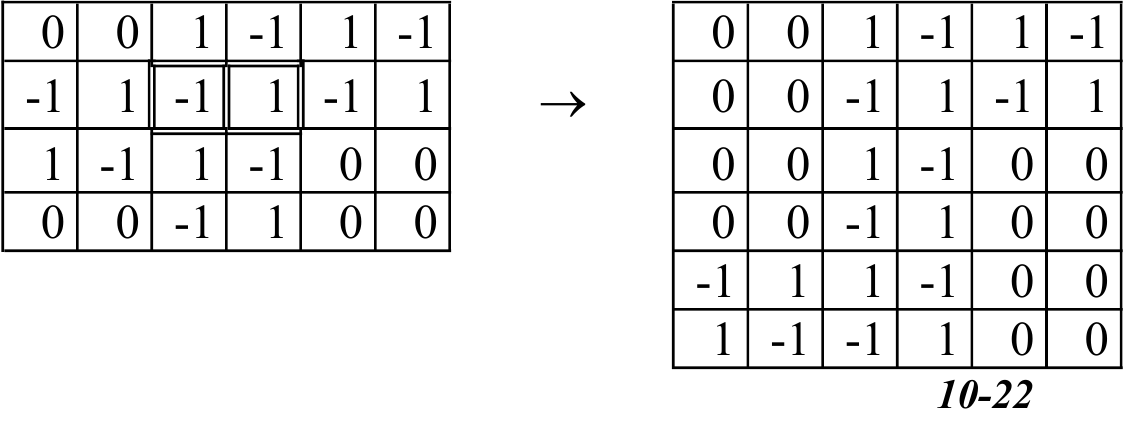}

The initial vectors. Corresponding graph and examples of sweeps.

\begin{figure}[ht] 
\center{\includegraphics[height=3.5cm]{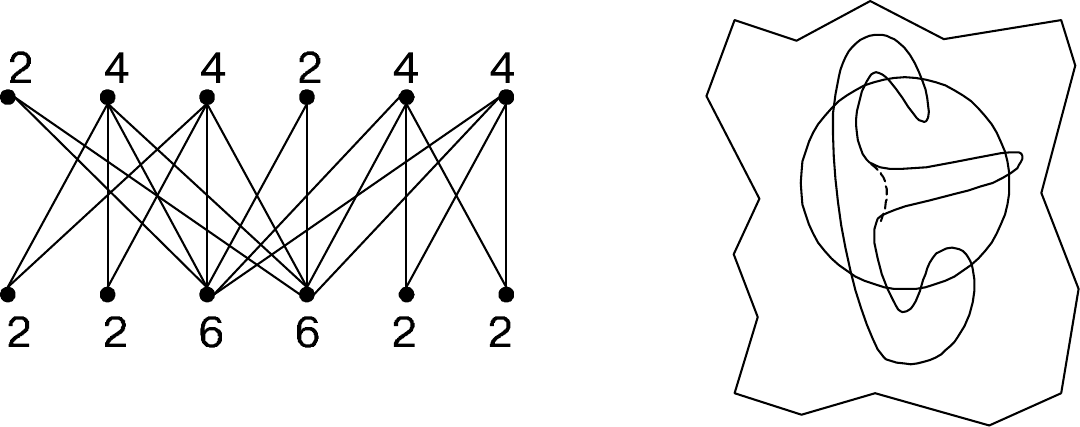}}
\label{p22}
\caption{ }
\end{figure} 

\subsection{ Now consider the intersection of the 8-5 matrix.}

     \subsubsection{ Case 1}
x -distribution is uniform.

\medskip                      
                       
\includegraphics[height=4.5cm]{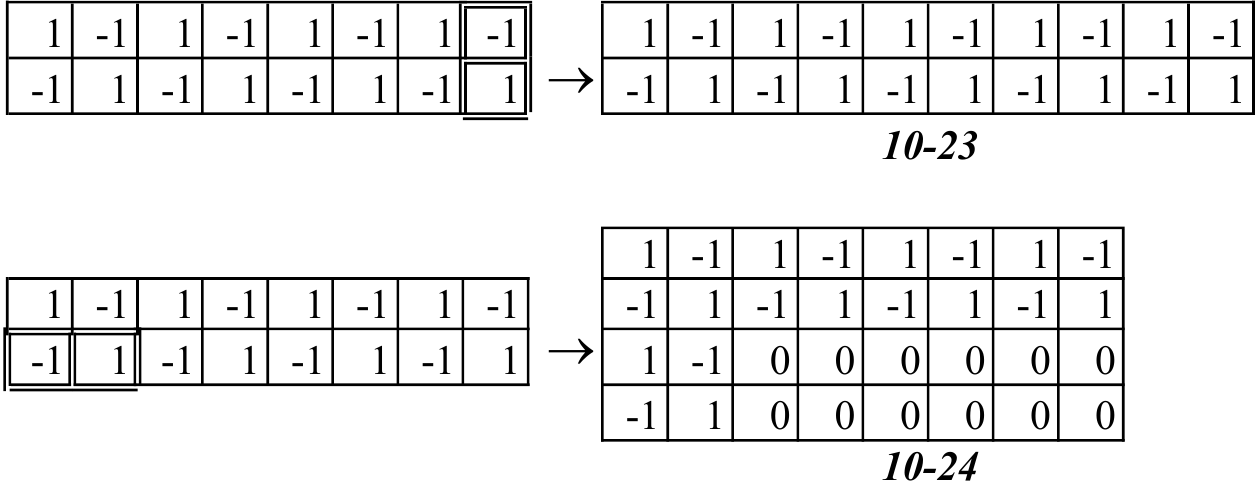}

The initial vectors 10-23. Corresponding graph and examples of sweeps.

\begin{figure}[ht] 
\center{\includegraphics[height=3.5cm]{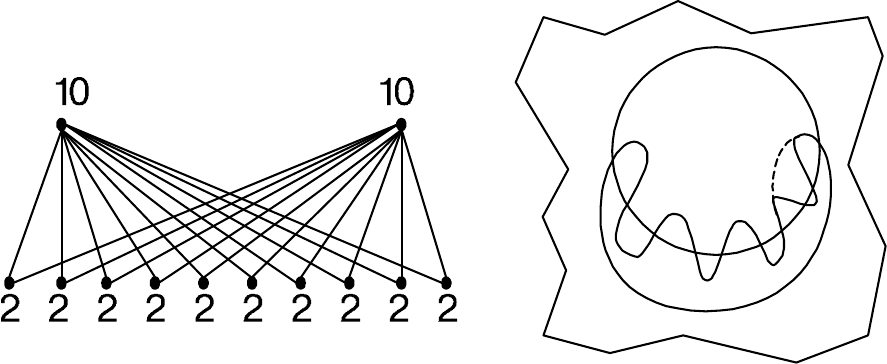}}
\label{p23}
\caption{ }
\end{figure}

    \subsubsection{ Case 2}
x -distribution is not uniform.

a)

\medskip

The initial vectors 10-24. Corresponding graph and examples of sweeps.

\begin{figure}[ht] 
\center{\includegraphics[height=3.5cm]{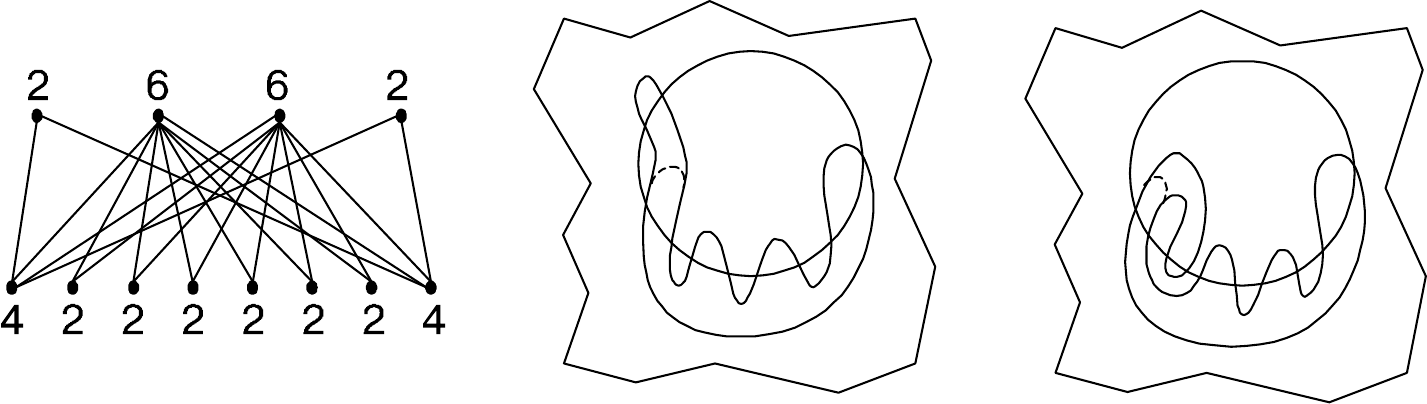}}
\caption{ }
\end{figure}

                       b)
\medskip                      
                       
\includegraphics[height=2.5cm]{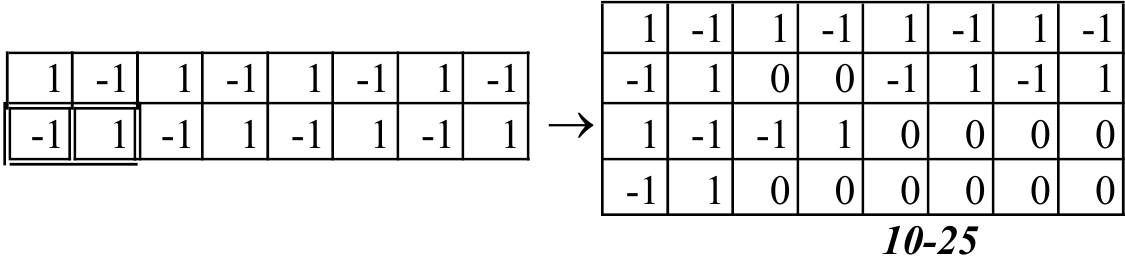}

The initial vectors. Corresponding graph and examples of sweeps.

\begin{figure}[ht] 
\center{\includegraphics[height=3.5cm]{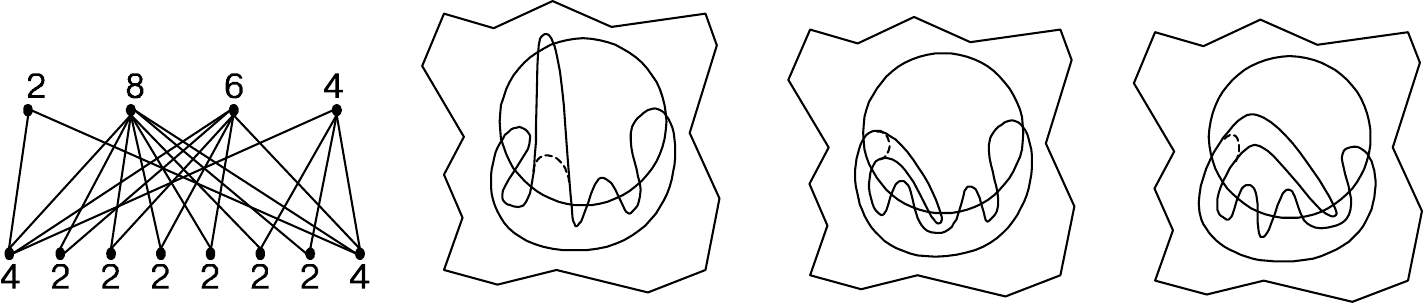}}
\caption{ }
\end{figure}

         \subsection{      Calculation of results}

If the defining vectors of the configurations differ, then the configurations also differ. Therefore, among the configurations given above, we will look for equivalent ones only among those that have the same defining vectors.

\begin{figure}[ht] 
\center{\includegraphics[height=5.0cm]{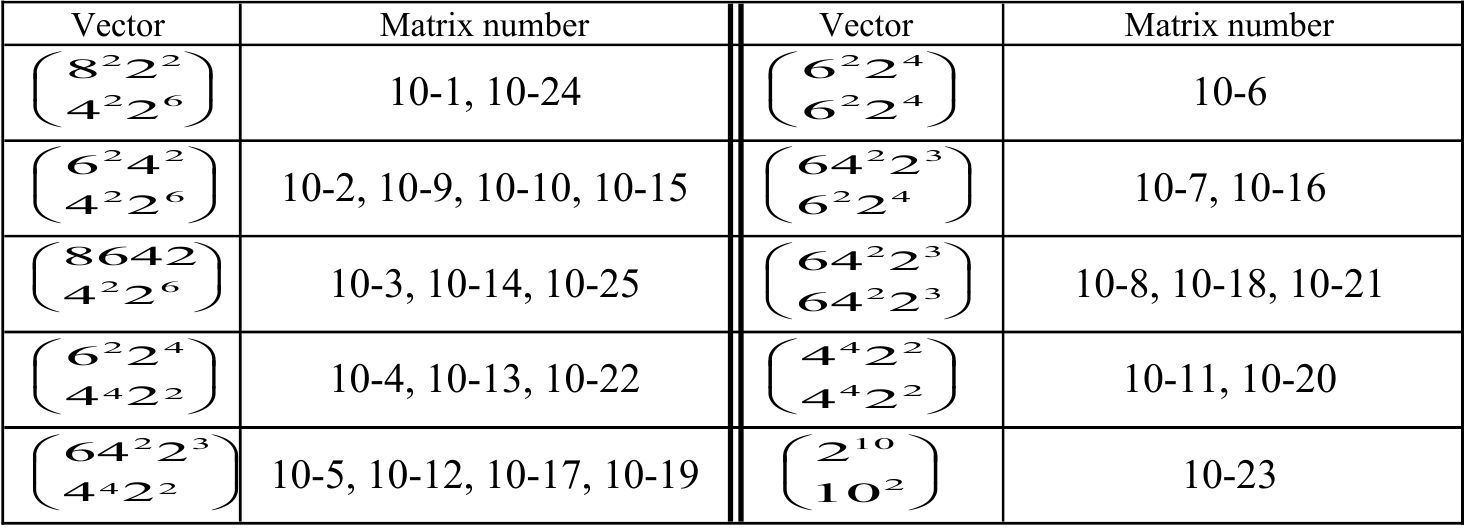}}
\end{figure}

Upon direct verification (one can equate to the non-equivalence of a graph or matrix, or to equate to the equivalence of a sweep) we have: 
matrices 10-1 and 10-24 are equivalent; 10-2 and 10-9; 10-10 and 10-15; 10-3, 10-14 and 10-25; 10-4 and 10-13; 10-5, 10-12, 10-17 and 10-19; 10-7 and 10-16;  10-18 and 10-21.
Now we can write out all configurations from 10 points: 10-1 (fig. 17), 10-2 (fig. 18),10-3 (fig. 19),10-4 (fig. 20),10-5 (fig. 21),10-6 (fig. 22),10-7 (fig. 23),10-8 (fig. 24),10-11 (fig. 27),10-18 (fig. 34),10-20 (fig. 36),10-22 (fig. 38),10-23 (fig. 39).

Thus, the following theorem holds.
\begin{theorem}
There are 13 configurations with 10 intersections: $N(10) =13$.
\end{theorem}

\section{Morse flow and the intersection of two closed curves on the sphere}

We consider Morse flow ( gradient-like Morse-Smale vector fields) on a 3-dimensional sphere, which have one singular point of index 1 and 2 and two sinks and sources each.

\begin{figure}[ht] 
\center{\includegraphics[height=3.5cm]{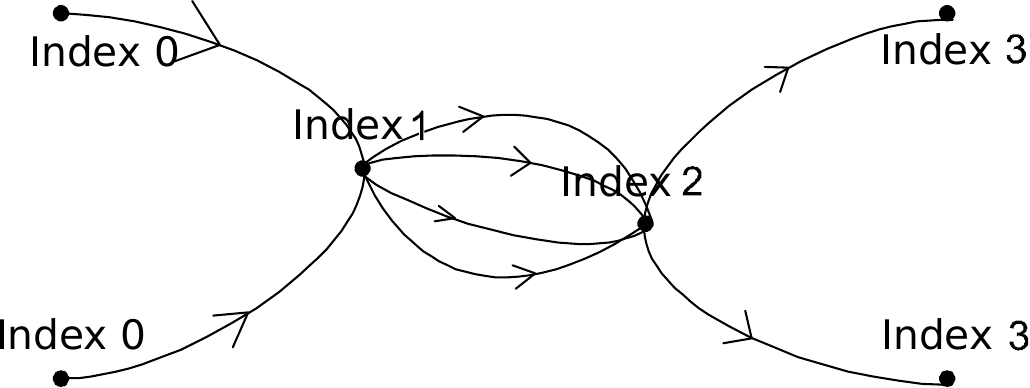}}
\caption{ }
\end{figure}

One trajectory from each of the sources enters into the point of index 1. From the point of index 2, one trajectory leads to each of the sources. Let $F_1$ be the regular neighborhood circle of the first two trajectories. Its boundary is topologically equivalent to $S^2$. Let $F_2$ be the regular neighborhood of their last two trajectories. Its boundary is also topologically equivalent to $S^2$.  Intersections of  the trajectories starting from the point of index 1 with the boundary of $F_1$ forms a closed curve. In the same way, a closed curve is formed on the boundary of $F_2$. We have: from each point of $\partial F_1$ a certain trajectory is going out, to each point of $\partial F_2$ some trajectory is included .

\begin{figure}[ht] 
\center{\includegraphics[height=4.5cm]{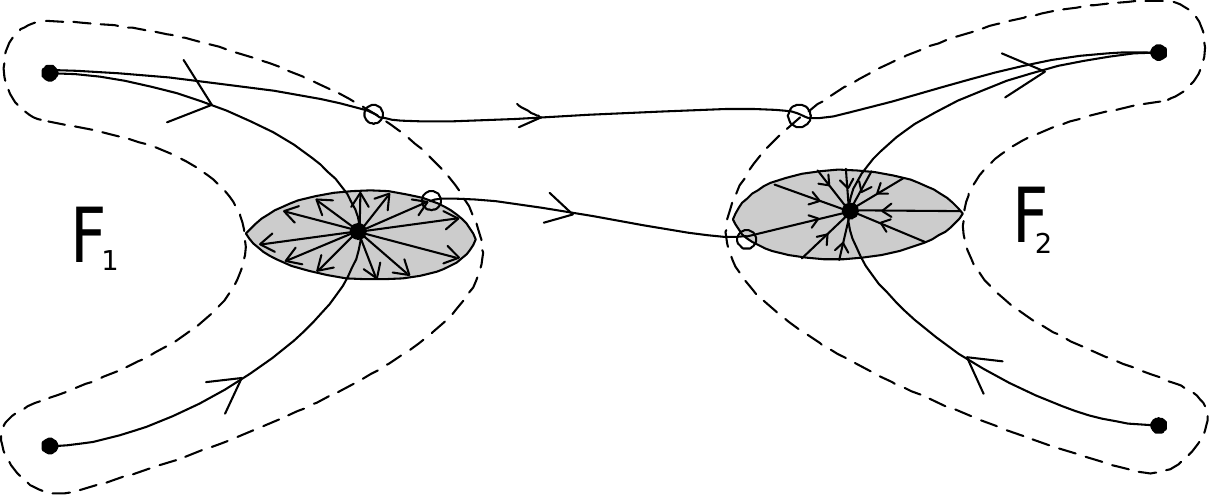}}
\caption{ }
\end{figure}

Thus, we can identify $\partial F_1$ and $\partial F_ 2$. Two closed curves corresponding to points with indices 1 and 2 remain on the formed sphere. The intersections of these curves were considered in the previous  sections of the paper. Each point of intersection of two curves corresponds to a saddle connection, trajectory connecting fixed points of index 1 and index 2. According to \cite{prish1998vek, prishlyak1999equivalence, prish2001top, prishlyak2002morse1,  Prishlyak2002,  prishlyak2002topological, prishlyak2005complete}, the configuration of two closed curves on the sphere defines the field with accuracy up to topological equivalence.

When classifying the intersections of closed curves on the sphere, we considered the configurations to be equivalent even if there was a homeomorphism that translates the first curve into the second and vice versa. When applied to flow classification, such configurations must be distinguished.

Therefore, the results of the previous section can be adjusted as follows: when calculating the number of topologically non-equivalent flows, each configuration from the first section is checked for symmetry with respect to the permutation of the curves. If the configuration turns out to be symmetric, then (by construction) one flow corresponds to it, otherwise, an asymmetric configuration corresponds to two topologically non-equivalent flows. Checking for symmetry of a configuration is checking the reduced matrix of the adjacency graph of this configuration. If the matrix remains equal to itself when multiplied by -1 (exchange of curves by places) with accuracy up to arbitrary permutation of rows and columns and rotation by $\pi /2$, then the matrix corresponds to a symmetric configuration, otherwise the matrix corresponds to an asymmetric configuration.

For cases with 2, 4, 6, and 8 intersection points (saddle connections), all configurations are symmetrical and the number of topologically non-equivalent fields is 1, 1, 2, 4, respectively. For the case with 10 intersection points, only the 10-18 matrix is asymmetric, which corresponds to two configurations:

\begin{figure}[ht] 
\center{\includegraphics[height=3.5cm]{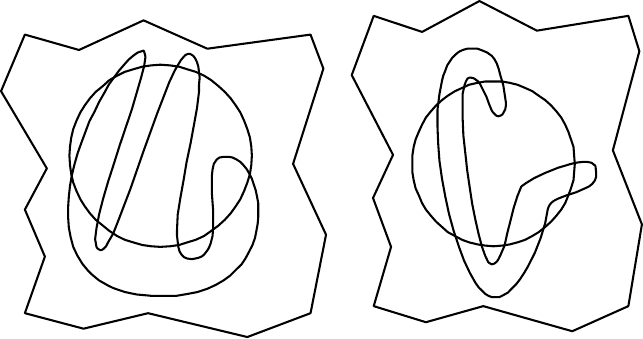}}
\caption{ }
\end{figure}

Thus, there are 14 topologically non-equivalent Morse flows with 10 saddle connections on $S^3$, in which there is one fixed point of index 1 and 2 and two sinks and sources each.

\section*{Conclusion} 

We found all possible configurations of two nested circles on a sphere with 2, 4, 6, 8 and 10 points of transversal intersection. The corresponding configuration numbers are: 1,1,2,4, 13. They specify the Morse flows on the three-dimensional sphere, respectively 1,1, 2, 4 and 14 flows. We hope that the obtained results can be generalized to other flows.

\bibliographystyle{plain}
\bibliography{prishd}

\textsc{Taras Shevchenko National University of Kyiv}

\textit{Email address:} \text{ bilun@knu.ua}   \ \ \ \ \ \ \ \ \
\textit{ Orcid ID:}  \text{0000-0003-2925-5392}

\textit{Email address:} \text{ prishlyak@knu.ua} \ \ \ \
\textit{ Orcid ID:} \text{0000-0002-7164-807X}

\end{document}